\documentclass[a4paper,fleqn,colorlinks,hypertexnames=false]{cas-sc}

\usepackage[authoryear,longnamesfirst]{natbib}
\usepackage{comment}
\usepackage{amsmath}
\numberwithin{equation}{section}

\def\tsc#1{\csdef{#1}{\textsc{\lowercase{#1}}\xspace}}
\tsc{WGM}
\tsc{QE}
\DeclareMathOperator{\sech}{sech}

\usepackage{float}
\usepackage{hyperref}

\begin{document}
\let\WriteBookmarks\relax
\def\floatpagepagefraction{1}
\def\textpagefraction{.001}

% Short title
\shorttitle{Multiple scales analysis of a nonlinear timestepping instability}    

% Short author
\shortauthors{B. A. Hyatt et-al.}  

% Main title of the paper
\title [mode = title]{Multiple scales analysis of a nonlinear timestepping instability in simulations of solitons}  

\author[1,2]{Benjamin A. Hyatt}[orcid=0000-0002-7401-3256]

% Corresponding author indication
\cormark[1]

% Email id of the first author
\ead{benjamin.hyatt@u.northwestern.edu}

% Address/affiliation
\affiliation[1]{organization={Department of Engineering Sciences and Applied Mathematics, Northwestern University},
            city={Evanston},
            state={IL},
            postcode={60208}, 
            country={USA}}

\affiliation[2]{organization={Center for Interdisciplinary Exploration and Research in Astrophysics, Northwestern University},
            city={Evanston},
            state={IL},
            postcode={60201}, 
            country={USA}}

\author[1,2]{Daniel Lecoanet}[orcid=0000-0002-7635-9728]

\author[2,3]{Evan H. Anders}[orcid=0000-0002-3433-4733]

\affiliation[3]{organization={Kavli Institute for Theoretical Physics, University of California, Santa Barbara},
            city={Santa Barbara},
            state={CA},
            postcode={93106},
            country={USA}}

\author[4,5]{Keaton J. Burns}[orcid=0000-0003-4761-4766]

\affiliation[4]{organization={Department of Mathematics, Massachusetts Institute of Technology},
            city={Cambridge},
            state={MA},
            postcode={02139},
            country={USA}}

\affiliation[5]{organization={Center for Computational Astrophysics, Flatiron Institute}, 
            city={New York},
            state={NY},
            postcode={10010},
            country={USA}}

% Footnote of the second author
%\fnmark[2]

% Email id of the second author
%\ead{}

% URL of the second author
%\ead[url]{}

% Credit authorship
%\credit{}

% Address/affiliation
%\affiliation[2]{organization={},
%            addressline={}, 
%            city={},
%          citysep={}, % Uncomment if no comma needed between city and postcode
%            postcode={}, 
%            state={},
%            country={}}

% Corresponding author text
\cortext[1]{Corresponding author}

% Footnote text
%\fntext[1]{}

% For a title note without a number/mark
%\nonumnote{}

\begin{abstract}
The susceptibility of timestepping algorithms to numerical instabilities is an important consideration when simulating partial differential equations (PDEs). Here we identify and analyze a pernicious numerical instability arising in pseudospectral simulations of nonlinear wave propagation resulting in finite-time blow-up. The blow-up time scale is independent of the spatial resolution and spectral basis but sensitive to the timestepping scheme and the timestep size. The instability appears in multi-step and multi-stage implicit-explicit (IMEX) timestepping schemes of different orders of accuracy and has been found to manifest in simulations of soliton solutions of the Korteweg-de Vries (KdV) equation and traveling wave solutions of a nonlinear generalized Klein-Gordon equation. Focusing on the case of KdV solitons, we show that modal predictions from linear stability theory are unable to explain the instability because the spurious growth from linear dispersion is small and nonlinear sources of error growth converge too slowly in the limit of small timestep size. We then develop a novel multi-scale asymptotic framework that captures the slow, nonlinear accumulation of timestepping errors. The framework allows the solution to vary with respect to multiple time scales related to the timestep size and thus recovers the instability as a function of a slow time scale dictated by the order of accuracy of the timestepping scheme. We show that this approach correctly describes our simulations of solitons by making accurate predictions of the blow-up time scale and transient features of the instability. Our work demonstrates that studies of long-time simulations of nonlinear waves should exercise caution when validating their timestepping schemes. 

\end{abstract}

% Use if graphical abstract is present
%\begin{graphicalabstract}
%\includegraphics{}
%\end{graphicalabstract}

% Research highlights
%\begin{highlights}
%\item 
%\item 
%\item 
%\end{highlights}

% Keywords
% Each keyword is seperated by \sep

\maketitle

% Main text

\section{Introduction}\label{sec:intro}

Distinguishing physical from computationally spurious behaviors is essential for interpreting simulations of partial differential equations (PDEs). One common source of spurious behavior is numerical instability. Avoiding numerical instability may place stringent restrictions on spatial and temporal resolutions, limiting the efficiency of numerical schemes. We therefore study numerical instabilities to learn their driving mechanisms and better discern which schemes are compatible with certain classes of problems. 

When a particular timestepping scheme is expected to cause an instability, it is standard practice to make a diagnosis using linear stability theory. A common example is performing a von Neumann stability analysis, which follows a strategy set forth in \cite{Charney1950} and formalized in many subsequent works (e.g., \citealt{Chan1984}, \citealt{Trefethen1994}, \citealt{Kumari2021}) to see how the temporal discretization evolves the Fourier modes of a perturbation to the numerical solution. While this approach is commonplace, it is not without limitations. For instance, when a timestepping scheme involves the evolution of non-normal linear operators, 
a von Neumann analysis may be insufficient to determine credible stability bounds without additional constraints on the operator (pseudo)spectrum (e.g., as discussed in \citealt{Scobelev1998}, \citealt{Reddy1990}). Furthermore, in nonlinear problems, performing a linear stability analysis requires deciding on a suitable linearization procedure, such as performing a diagonalization by way of freezing non-constant coefficients (e.g., \citealt{Fornberg1978}, \citealt{Fuhrman2004}). Hence, even in problems of modest complexity, basic linear stability analyses run the risk of neglecting or overlooking aspects of the discretization which play a critical role in the numerical instability. 

A complementary perspective is that numerical instabilities can develop from the violation of CFL-type conditions. First recognized in \cite{Courant1928}, CFL conditions are restrictions for stability on the timestep size, $\Delta t$, which ensure the physical solution propagates no faster than the speed at which the numerical scheme can send information throughout the domain. This accounts for the CFL heuristic that the stable timestep size in explicit timestepping schemes is limited by the smallest spatial scales in the simulation. On the other hand, when using implicit timestepping schemes (or, for instance, when using spatial discretizations that share information globally, e.g., Galerkin spectral methods), discrete information is shared everywhere on each step, so this statement of CFL is inapplicable. However, CFL may be used more loosely to refer to any stability restriction on $\Delta t$ that can be expressed in terms of the spatial resolution and an advection speed. Such CFL-type conditions are frequently established by way of a von Neumann analysis and other applications of linear theory (e.g., \citealt{Gottlieb1991}, \citealt{Grooms2011}, \citealt{Deriaz2020}). 

While the standard tools of linear stability theory are widely applicable, they may come up short in the description of numerical instabilities in nonlinear simulations. The impact of nonlinearities on numerical instability is perhaps best understood in terms of the spatial discretization, as there is a class of instabilities known to arise from the aliasing of spatial modes in the evaluation of nonlinear terms. Much work in this area stems from \cite{Phillips1959}, which articulates how schemes can mistakenly allow nonlinear interactions between small-scale modes to affect larger scales in the flow, and thereby create feedback loops that lead to blow-up. Subsequent works demonstrate how this aliasing instability primarily depends on the initial noise level and the smallest resolved length scales, which together determine a slow time scale over which energy accumulates at small scales before exceeding a critical amplitude and quickly blowing up thereafter (e.g., \citealt{Newell1977}, \citealt{Briggsetal1983}, \citealt{Mitchell1988}). Furthermore, as demonstrated by \citealt{Newell1977}, restriction of the timestep size according to a linear analysis of the temporal discretization is usually incapable of preventing the instability. Most remedies---such as smoothing, dealiasing, use of energy-conserving spatial discretizations, and so on---instead focus on suppressing or altogether eliminating aliasing errors to thereby avoid their instabilities (e.g., Chapter 11 in \citealt{Boyd2000}, \citealt{Bardos2015}). Hence, this type of instability sensitively depends on the numerical strategies used to contend with small spatial scales in the flow, whereas the choice of timestepping scheme is viewed as somewhat ancillary. 

In this work, we identify a nonlinear numerical instability driven exclusively by the temporal discretization---i.e., the spatial discretization is inconsequential---in nonlinear simulations of wave propagation. The instability is highly sensitive to the choice of timestep size in a way that depends on the particular timestepping scheme and its accumulation error. Rather than exciting energy at particular spatial scales, this instability is marked by smooth changes in the solution envelope. Beginning at early times, these changes result in the numerical and exact solutions gradually moving out of phase, accompanied by deviations in the envelope amplitude and length scale. 
Moreover, \textit{linear theory fails} to predict these developments, as at the earliest times in our simulations and in the limit of small timestep size, linear sources of spurious growth are overpowered by nonlinear contributions to the error. The progression of the spurious error growth accelerates over time, and in many cases culminates in finite-time blow-up.  

This nonlinear timestepping instability may play an important role in a wide class of problems involving nonlinear waves. In particular, we have observed this instability in the long-time evolution of soliton solutions of the Korteweg-de Vries equation~\eqref{eq:KdV} and traveling wave solutions of a generalized Klein-Gordon equation with a cubic nonlinearity (e.g., \citealt{Liu2001}). We hypothesize that the instability may emerge in other similar contexts, e.g., in Boussinesq equations (\citealt{Liu2001}) or Kuramoto–Sivashinsky equations (\citealt{Nickel2007}), which likewise admit exact soliton and traveling wave solutions. Propagating solitons and other waveforms arise in countless areas of mathematical physics, often with an interest in simulating their dynamics over long time scales. For example, the seminal work of \cite{ZabuskyKruskal1965} studies a ``recurrence'' phenomenon where an initial excitation of solitons return nearly to their original state after a long time has passed. Other studies focus on the emergence of new features over long time scales, such as changes in soliton parameters when the flow is forced by random perturbations (e.g., \citealt{Garnier2001}), the generation of soliton ensembles from periodic initial data (e.g., \citealt{Salupere2003}), the long-term effects of dissipation on soliton dynamics (e.g., \citealt{Christov1995}), and the aftermath of multiphase wave interactions (e.g., \citealt{Mao2023}). Since these lines of inquiry rely heavily on numerical experimentation to identify new behaviors and verify analytic predictions playing out over many iterations, it is important that their timestepping schemes avoid insidious slow-growing instabilities.

In this paper, we use a multi-scale asymptotic analysis to describe this nonlinear timestepping instability. Notably, while the numerical solution exhibits near-analytic behavior on short time scales, there is clear evolution of the numerical error with respect to a long time scale whose duration is found to be a function of the timestep size. This separation of time scales motivates us to perform an asymptotic analysis involving {\it multiple time scales}. The asymptotic analysis also crucially retains nonlinear features of the problem, essential for this instability which we find cannot be reproduced by a linear analysis.

The introduction of multiple time scales is a well-known technique for approximating the solution to nonlinear problems, and is useful for finding approximations which converge uniformly in time (e.g., Chapter 3 of \cite{Holmes2013}). While the method is commonly applied to continuous-time problems, the present work leverages a multiple scales analysis in the discrete-time context of timestepping. We believe our particular use of this method in the description of a timestepping instability is novel. However, we highlight that the aforementioned work of \cite{Newell1977} in fact used a similar approach to study the aliasing instability, and works such as \cite{Hoppensteadt1977} have briefly motivated the analysis of numerical schemes as a potential application of the multiple scales strategy. 
 
This study demonstrates that the observed timestepping instability cannot be explained by conventional and heuristic methods from linear theory, but is rather a consequence of the nonlinear behavior of the temporal discretization. Section~\ref{sec:sol} introduces simulations of solitons---which serve as a case study to explore the instability---and shows the instability is not related to the traditional CFL criterion. Section~\ref{sec:lin-stab} then performs a von Neumann-type linear stability analysis and compares the predicted modal behavior against our simulations. Finding that the linear analysis does not describe the observed timestepping instability, we motivate an analysis which captures the nonlinear behavior of the numerical solution. Section~\ref{sec:ms} develops a multiple scales analysis capable of predicting the time scale of the instability and its spurious features. We compare the predictions of the analysis to simulations in Section~\ref{sec:results}, and offer closing remarks in Section~\ref{sec:conc}.

\section{Timestepping instability in simulations of solitons}\label{sec:sol}

Nonlinear waves, including solitons, appear in many distinct physical systems. One important model which admits solitons is the Korteweg-de Vries (KdV) equation 
\begin{equation} \label{eq:KdV}
    \partial_t u + \alpha \partial_x^3 u = - u \partial_x u.%,
\end{equation} 
Introduced in \cite{KortewegdeVries1895}, \eqref{eq:KdV} is well known to have exact solutions taking the form of a solitary traveling wave, i.e., a soliton. While we have observed similar timestepping instabilities in multiple equations admitting traveling wave solutions, we will focus here on an analysis of the KdV problem. The solutions of~\eqref{eq:KdV} which, at any fixed time $t$, decay at infinity are solitons of the form
\begin{equation} \label{eq:soliton}
    u(x,t) = 3c\sech\left(\sqrt{\frac{c}{4\alpha}}\left(x-ct\right)\right)^2
\end{equation}
where the constant $c > 0$ sets the amplitude and speed of propagation. The soliton length scale is further set by $c$ and the dispersion parameter $\alpha > 0$. As seen in \eqref{eq:soliton}, smaller values of $\alpha$ correspond to solitons with smaller length scales. It is convenient to note that the full width at half maximum of the soliton in \eqref{eq:soliton} is given by \begin{equation}\label{eq:fwhm}
    x_{\rm FWHM} = 2\sqrt{\frac{4\alpha}{c}}\ln\left(1+\sqrt{2}\right).
\end{equation}
There are more complex varieties of nonlinear solutions which exactly solve KdV-type equations, such as $N$-soliton solutions and their generalizations (e.g., \cite{Deift1979}, \cite{Gesztesy1992}), but we focus on the simple case of the single soliton here.

This work is primarily interested in the numerical behavior of solitons over long time scales. To study this long-term behavior, \eqref{eq:KdV} is evolved over the periodic domain $x \in [-L/2, L/2]$, taking $L = 10$. Solitons are initialized at $t = 0$ at the center of the domain according to \eqref{eq:soliton}, and thereafter they propagate in the $+x$ direction with continuous passage through the periodic boundaries. While \eqref{eq:soliton} is not exactly periodic, it is approximately periodic provided that the length scale set by $c$ and $\alpha$ is sufficiently small. We take $c = 0.5$ such that the expected time scale for one domain crossing in each simulation is $T_c = L/c = 20$. To ensure our soliton is both approximately periodic and spatially resolved, we consider values of $\alpha$ in the range of $0.002$ to $0.03$. This corresponds to $x_{\rm FWHM}$ less than one tenth of the domain size $L$, and initial values of $u \lesssim 10^{-12}$ at the boundaries $x = \pm L/2$.

We evolve our simulations of solitons with pseudospectral methods using the open-source pseudospectral code {\tt Dedalus} (\citealt{Burns2020}). Timesteps are taken using mixed implicit-explicit (IMEX) schemes where the nonlinear advection term in \eqref{eq:KdV} is evaluated explicitly and the linear dispersion term is evaluated implicitly. We employ $3/2$-dealiasing to prevent the quadratic nonlinearity in~\eqref{eq:KdV} from producing aliasing errors. There are several studies in the literature which consider the numerical stability of IMEX schemes paired with spectral methods (e.g., \cite{Ascher1995}), including in the context of nonlinear problems with dispersion (e.g., \citealt{Fornberg1978}, \citealt{Fornberg1999}, and \citealt{Grooms2011}). Each of the references given here uses linear theory to assess numerical stability, which to the best of our knowledge is commonplace. 

We consider IMEX schemes of several different orders of accuracy drawn from two classes, including semi-implicit backward differencing (SBDF) schemes and mixed explicit, diagonally-implicit (ERK+DIRK) Runge-Kutta schemes. The specific schemes are provided in Appendix~\ref{sec:IMEX} and are used in the analysis in both Section~\ref{sec:lin-stab} and Section~\ref{sec:ms}.

We use two different spatial discretizations for this spatially periodic problem; $u$ is expanded either in terms of a real Fourier basis of sines and cosines, or a basis of Chebyshev $T$ polynomials. In the latter case, periodicity is enforced via the boundary conditions
\begin{equation}\label{eq:chebyshevBCs}
    u{\left(-\frac{L}{2}, t\right)} = u{\left(\frac{L}{2}, t\right)}, \quad \partial_x u{\left(-\frac{L}{2}, t\right)} = \partial_x u{\left(\frac{L}{2}, t\right)}, \quad \partial_x^2 u{\left(-\frac{L}{2}, t\right)} = \partial_x^2 u{\left(\frac{L}{2}, t\right)}
\end{equation}
using a generalized tau method as described in~\cite{Burns2024}. Although using Fourier modes is standard, we also use Chebyshev polynomials to study the stability properties of the two spatial discretizations. 

\begin{figure}[t]
    \centering
    \includegraphics{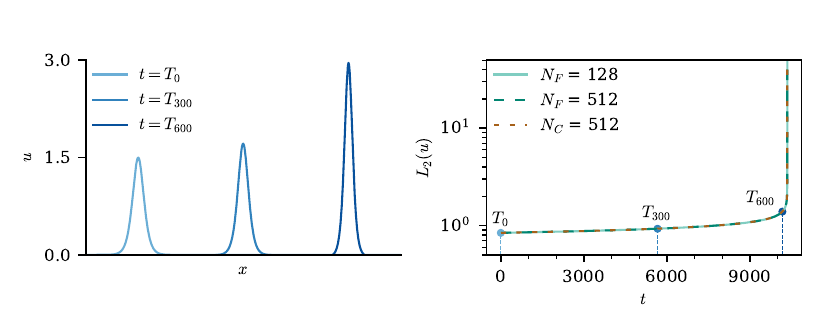}
    \caption{\textbf{Evolution of the nonlinear timestepping instability.} The numerical solution experiences unstable growth in the $L_2$ norm while maintaining the general shape of a soliton. The instability consists of the soliton increasing in amplitude, shrinking in width, and gaining speed. The left panel depicts three snapshots  %in time 
    of the soliton when $\alpha = 0.00697$ in a simulation with $\Delta t = 0.00493$ using the SBDF2 timestepping scheme and a Fourier basis with $N = 512$.
    Snapshots are recorded at $T_0 = 0$, $T_{300} = 5657$, %$T_{300} = 5656.914$, 
    and $T_{600} = 10160$, %$T_{600} = 10164.037$, 
    which denote the times at which the soliton has traversed the domain 0, 300, and 600 times, respectively. The right panel tracks the evolution of the $L_2$ norm of the numerical solution for three different spatial discretizations---including Fourier and Chebyshev bases with different numbers of modes---and shows that they are indistinguishable.
    }
    \label{fig:1-inst-vis}
\end{figure}

We find that our simulations of solitons are unstable for most choices of IMEX timestepping scheme. Figure~\ref{fig:1-inst-vis} provides a visual depiction of the evolution of this numerical instability, here produced by the second-order accurate semi-implicit backward difference scheme, SBDF2. The soliton first propagates in the $+x$ direction, passing through the periodic boundaries hundreds of times. Snapshots of the numerical solution reveal that the soliton gradually deviates from the exact solution by increasing in amplitude and decreasing in width, but overall retaining a shape similar to the initial soliton data. The increase in amplitude is accompanied by a proportional increase in the propagation speed over time, resulting in the numerical solution growing out of phase with the exact solution. All the while, the $L_2$ norm of the numerical solution grows over time. The growth in the $L_2$ norm accelerates, eventually leading to blow-up.

Notably, the instability tracked in Figure~\ref{fig:1-inst-vis} follows a similar evolution for \textit{different spatial resolutions}.
Changing the spatial resolution or changing the basis functions (e.g., Fourier or Chebyshev) does not alter the development of the instability.
This rules out an instability of CFL type, as the instability is independent of the spatial discretization.
We find ubiquitous unstable behavior that is sensitive to the timestep size, $\Delta t$, and the dispersion parameter, $\alpha$. The dependence on $\Delta t$ indicates this instability is related to the temporal discretization, i.e., this is a timestepping instability. 

Moreover, Table~\ref{tbl:abbrev-survey} reports the blow-up times as a function of spatial discretization, $\alpha$, and $\Delta t$ in simulations using the SBDF2 timestepping scheme. As the number of basis functions $N$ is increased, $T_{\rm blowup}$ appears to converge, suggesting this is an instability of the semi-discrete in time problem.

On the other hand, when one of $\Delta t$ or $\alpha$ is varied there is a clear and significant change in $T_{\rm blowup}$. In particular, smaller values of the timestep size $\Delta t$ and larger values of the dispersion parameter $\alpha$ correspond to a slower numerical instability and later $T_{\rm blowup}$. Thus, reducing the fixed timestep size slows down {\it but does not prevent} this numerical instability.

\begin{table}%[h]
\caption{\textbf{Measured blow-up times in simulations of solitons using the SBDF2 timestepping scheme}. Here we report blow-up times $T_{\rm blowup}$ for various spatial discretizations by changing the spectral basis and number of basis functions. We surveyed the blow-up times for a range of logarithmically spaced values of the dispersion parameter $\alpha$ between $0.002$ and $0.03$, and a range of logarithmically spaced fixed timestep sizes $\Delta t$ between $0.0006$ and $0.05$. This Table reports on a small subset of these results, which are representative of the behavior over the full parameter space. We find that the blow-up time is insensitive to the choice of spectral basis or number of basis functions $N$. There is strong dependence, however, on the dispersion parameter $\alpha$ and the timestep size $\Delta t$.}\label{tbl:abbrev-survey}
\begin{tabular*}{\tblwidth}{@{} L | L L L L @{}}
\toprule
 Basis & $N$ & $\alpha$ & $\Delta t$ & $T_{\rm blowup}$ \\ % Table header row
\midrule \midrule
Fourier & $N_F = 128$ & $0.00697$ & $0.00324$ & $36779.374$ \\
Chebyshev & $N_C = 128$ & $0.00697$ & $0.00324$ & $37060.177$\\
Fourier & $N_F = 256$ & $0.00697$ & $0.00324$ & $36771.441$ \\
Chebyshev & $N_C = 256$ & $0.00697$ & $0.00324$ & $36771.645$ \\
Fourier & $N_F = 512$ & $0.00697$ & $0.00324$ & $36771.441$\\
Chebyshev & $N_C = 512$ & $0.00697$ & $0.00324$ & $36771.441$ \\
\midrule \midrule
Fourier & $N_F = 128$ & $0.00697$ & $0.00609$ & $5472.393$\\
Chebyshev & $N_C = 128$ & $0.00697$ & $0.00609$ & $5517.486$\\
Fourier & $N_F = 256$ & $0.00697$ & $0.00609$ & $5471.084$\\
Chebyshev & $N_C = 256$ & $0.00697$ & $0.00609$ & $5471.127$ \\
Fourier & $N_F = 512$ & $0.00697$ & $0.00609$ & $5471.084$ \\
Chebyshev & $N_C = 512$ &$0.00697$ & $0.00609$ & $5471.084$ \\
\midrule \midrule
Fourier & $N_F = 128$ & $0.01312$ & $0.00324$ & $129853.382$\\
Chebyshev & $N_C = 128$ & $0.01312$ & $0.00324$ & $129896.992$\\
Fourier & $N_F = 256$ & $0.01312$ & $0.00324$ & $129852.638$\\
Chebyshev & $N_C = 256$ & $0.01312$ & $0.00324$ & $129852.655$\\
Fourier & $N_F = 512$ & $0.01312$ & $0.00324$ & $129852.638$\\
Chebyshev & $N_C = 512$ & $0.01312$ & $0.00324$ & $129852.642$\\
\bottomrule
\end{tabular*}
\end{table}

\section{Linear stability}\label{sec:lin-stab}

\subsection{Linear analysis of the KdV equation}\label{sec:lin-stab-setup}

At early times in our simulations, the numerical solution is close to the exact soliton, $u_{\rm ex}$, given in~\eqref{eq:soliton}. Writing $u = u_{\rm ex} + v$, we can model the deviation from the exact soliton $v$ with the linearized equation
\begin{equation}\label{eq:KdV-lin}
\partial_t v + \alpha \partial_x^3 v = -u_{\rm ex} \partial_x v - \left(\partial_x u_{\rm ex} \right)v, \quad -\frac{L}{2} \le x \le \frac{L}{2},\quad t \ge 0,
\end{equation}
where the nonlinear term $-v\partial_x v$ has been assumed to be negligibly small and therefore dropped. Later we will show that this assumption may be invalid depending on the choice of timestepping scheme, even in the limit as $\Delta t \rightarrow 0$.

An alternative choice to~\eqref{eq:KdV-lin} would be to take the frozen coefficients approach and linearize the KdV equation about a constant background solution. In that case the stability analysis is related to the standard linear stability regions of each timestepping scheme. We address this method in Appendix~\ref{sec:regions} but find that it offers limited insights when compared to the current approach. 

Here we assess the stability of the perturbation $v$ by studying the evolution of its spatial discretization, given by
\begin{equation}\label{eq:sd}
    v(x, t) = \sum_{k = -N/2 + 1}^{k = N/2} \mathbf{v}_k(t) \exp{\left(i \frac{2\pi k}{L} x\right)} 
\end{equation}
where $\mathbf{v} \in \mathbb{C}^N$ is the vector of Fourier coefficients. %Here we assess the stability of the perturbation $v$ by assuming its evolution under the discretization of \eqref{eq:KdV-lin} can be expressed as 
We then assume the time evolution of each mode $\mathbf{v}_k(t)$ is exponential with a growth rate $\lambda$ to be determined and amplitude $\mathbf{v}_k$, such that the perturbation can be expressed as
\begin{equation}\label{eq:VN-ansatz}
    v(x, t_n) = \sum_{k = -N/2 + 1}^{k = N/2} \mathbf{v}_k \exp{(\lambda t_n)} \exp{{\left(i \frac{2\pi k}{L}(x - ct_n)\right)}} = \sum_{k = -N/2 + 1}^{k = N/2} \left(\mathbf{S}^n \mathbf{v}\right)_k  \exp{(\lambda t_n)}\exp{{\left(i \frac{2\pi k}{L} x\right)}}.
\end{equation}
Since the soliton $u_{\rm ex}$ translates by $c\Delta t$ to the right each timestep, we look for perturbations $v$ moving with the same speed. We capture this in~\eqref{eq:VN-ansatz} by defining a diagonal shift operator $\mathbf{S} = \text{diag}(\exp{(-ic\Delta t k)}) \in \mathbb{C}^{N\times N}$ acting on $\mathbf{v}$ each timestep. This lends to a clear interpretation of growth or decay of perturbations in the moving frame of the soliton.

Whereas the Fourier mode expansion represents the spatial discretization of the linearized dynamics, the time discretization is given by one of the IMEX timestepping schemes introduced in Appendix~\ref{sec:IMEX}. For each scheme we construct an eigenvalue problem whose eigenvalues $\sigma \equiv e^{\lambda \Delta t}$ represent amplification factors. The modulus $|\sigma|$ of each eigenvalue determines the change in the norm of a corresponding eigenmode in Fourier space, with $|\sigma| \le 1$ necessary for stability. 

For $s$-step SBDF schemes, we determine eigenvectors $\mathbf{w}^{(n)} \equiv \left(\mathbf{v}^{(n)}, \mathbf{v}^{(n+1)}, \dots, \mathbf{v}^{(n+s-1)} \right)^T \in \mathbb{C}^{sN}$  and their eigenvalues $\sigma$ by solving the $sN \times sN$ linear system
\begin{equation}\label{eq:sbdf-evp}
    \begin{aligned}
        %&\sigma \mathbf{S}\mathbf{v}^n - \mathbf{v}^{n+1} = 0, \\
        %&\sigma \mathbf{S}\mathbf{v}^{n+1} - \mathbf{v}^{n+2} = 0, \\ 
        %& \vdots \\
        &\sigma \mathbf{S}\mathbf{v}^{(n+i-1)} - \mathbf{v}^{(n+i)} = 0, \\
        & \vdots \quad (i = 1, \dots, s-1) \\
        &\sigma{\left(a_s\mathbf{I} + \Delta t \gamma_s \alpha \mathbf{D}^3 \right)}\mathbf{S}\mathbf{v}^{(n+s-1)} - \sum_{i=0}^{s-1}\biggl(-a_i\mathbf{I} - \Delta t \gamma_i \alpha \mathbf{D}^3 - \Delta t \beta_i \left(\left[u_{\rm ex}\right]\mathbf{D} + \left[\partial_x u_{\rm ex}\right]\right)\biggr)\mathbf{v}^{(n+i)} = 0,
    \end{aligned}
\end{equation} 
where $\mathbf{D} = \text{diag}(ik) \in \mathbb{C}^{N\times N}$ is the Fourier differentiation matrix, and $\mathbf{I}$ is the $N \times N$ identity matrix. Similarly, for $q$-stage ERK+DIRK schemes, we determine eigenvectors $\mathbf{w}^{(n)} \equiv \left(\mathbf{v}^{(n,0)}, \mathbf{v}^{(n,1)}, \dots, \mathbf{v}^{(n,q-1)}\right)^T \in \mathbb{C}^{qN}$ and their eigenvalues $\sigma$
by solving the $qN \times qN$ linear system
\begin{equation}\label{eq:rk-evp}
    \begin{aligned}
        &\left(\mathbf{I} + \Delta t \tilde{a}_{i, i} \alpha \mathbf{D}^3\right)\mathbf{v}^{(n,i)} - \mathbf{v}^{(n,0)} - \Delta t \sum_{j=0}^{i-1}\biggl(-\tilde{a}_{i,j}\alpha\mathbf{D}^3 -\hat{a}_{i,j}\left(\left[u_{\rm ex}(x,t_n + \tilde{c}_j \Delta t)\right]\mathbf{D} + \left[\partial_x u_{\rm ex}(x,t_n + \tilde{c}_j \Delta t)\right]\right)\biggr)\mathbf{v}^{(n,j)} = 0, \\
        & \vdots \quad (i = 1, \dots, q-1) \\
        &\sigma{\left(\mathbf{I} + \Delta t \tilde{a}_{q, q} \alpha \mathbf{D}^3\right)}\mathbf{S}\mathbf{v}^{(n,0)} - \mathbf{v}^{(n,0)} - \Delta t \sum_{j=0}^{q-1}\biggl(- \tilde{a}_{q,j}\alpha\mathbf{D}^3 -\hat{a}_{q,j}\left(\left[u_{\rm ex}(x,t_n + \tilde{c}_j \Delta t)\right]\mathbf{D} + \left[\partial_x u_{\rm ex}(x,t_n + \tilde{c}_j \Delta t)\right]\right)\biggr)\mathbf{v}^{(n,j)} = 0.
    \end{aligned}
\end{equation}
The eigenvalue problems differ from each other primarily in that $\sigma\mathbf{S}$ appears in each line of~\eqref{eq:sbdf-evp} due to the multi-step structure of the scheme, whereas $\sigma\mathbf{S}$ only appears once in~\eqref{eq:rk-evp} in the last line corresponding to the final stage. Equations~\eqref{eq:sbdf} and \eqref{eq:rk} in Appendix~\ref{sec:IMEX} respectively describe how $a_i$, $\beta_i$, and $\gamma_i$ specify the particular SBDF scheme and $\tilde{a}_{i,j}$ and $\tilde{c}_j$ specify the particular ERK+DIRK scheme. The particular values of these coefficients are stated in the Supplemental Material and can also be found in the references provided in Appendix~\ref{sec:IMEX}.

\subsection{Results of the linear stability analysis}\label{sec:lin-stab-results}

Linear theory predicts an SBDF or ERK+DIRK scheme to be unstable if there are eigenmodes with $|\sigma| > 1$. If we do find unstable modes, then we should verify their applicability by checking how their predicted growth rates compare to the corresponding numerical simulation. While the results of this linear analysis correspond to the numerical solution of the linearized initial value problem \eqref{eq:KdV-lin}, we will compare the predictions against the fully nonlinear initial value problem \eqref{eq:KdV} posed on the periodic finite interval and with soliton initial conditions. 

\begin{figure}[b]
    \centering
    \includegraphics{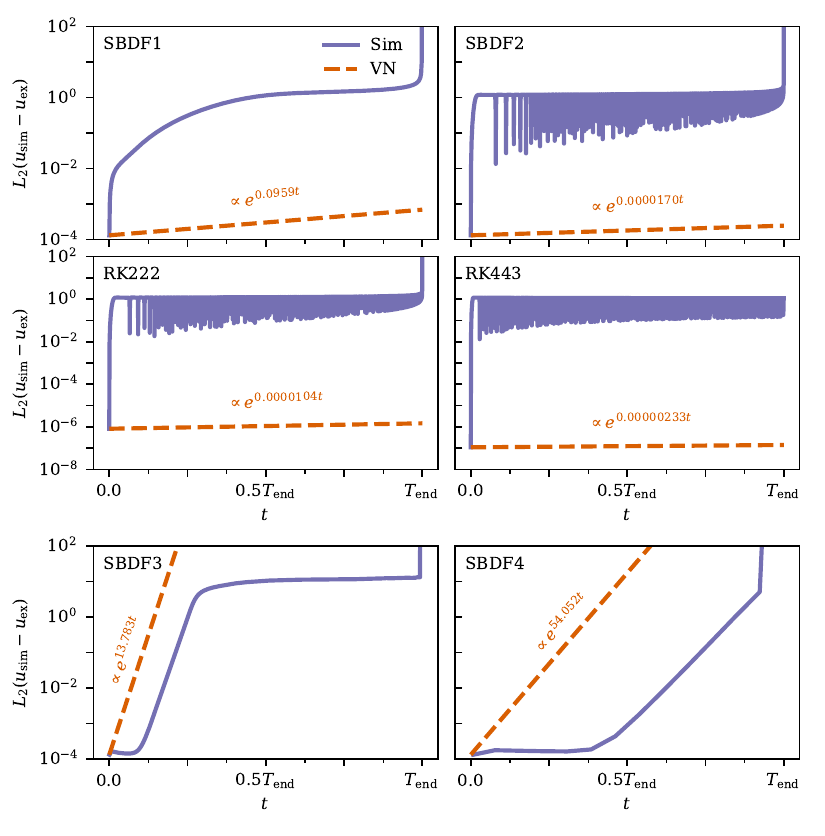}
    \caption{\textbf{Error in simulations (Sim) of solitons compared to predictions from von Neumann (VN) analysis}. We find linear instabilities in SBDF3 and SBDF4 in good agreement with the VN  predictions. By contrast, SBDF1, SBDF2, and RK222 exhibit fast non-modal growth at early times as the numerical solution quickly becomes out of phase with the exact solution. RK443 exhibits similar behavior but, unlike the preceding schemes, is stable and does not result in blow-up. Shown are the $L_2$ norms of the error 
    in pseudospectral simulations evolved using a timestep size of $\Delta t = 0.00324$ and dispersion parameter $\alpha = 0.00697$, alongside the norms of the fastest growing modes predicted by the von Neumann analysis. Each panel corresponds to a different timestepping scheme. Errors are plotted from $t = 0$ to $t = T_{\rm end}$ where $T_{\rm end}$ is the time at which the simulation terminated. The cause of termination is either due to blow-up (SBDF1: $T_{\rm end} = 17.35$, SBDF2: $T_{\rm end} = 36770$, SBDF3: $T_{\rm end} = 4.510$, SBDF4: $T_{\rm end} = 0.4367$, and RK222: $T_{\rm end} = 56680$) or reaching a user-specified time-limit (RK443: $T_{\rm end} = 3(30030/77069) \Delta t^{-3} \alpha^2 c_0^{-6} = 107400$). The errors in SBDF2, RK222, and RK443 momentarily ``dip'' whenever the numerical solution is in phase with the exact solution; since the numerical solution gains speed over time, it repeatedly laps the exact solution.
    }
    \label{fig:figure-2-evp-ivp}
\end{figure}

We solve the eigenvalue problems for SBDF1 (Forward-backward Euler), SBDF2, SBDF3, SBDF4, RK222, and RK443 for many fixed pairs of dispersion parameters $\alpha$ and timestep sizes $\Delta t$ sampled logarithmically from the ranges $0.003 \le \alpha \le 0.02$ and $0.0006 \le \Delta t \le 0.05$. Appendix~\ref{sec:evp-examples} provides example eigenspectra for each scheme and discusses our criteria for finding well-resolved, non-spurious eigenmodes. %discusses some general patterns we have observed. 

For all six schemes, we find that there are \textit{multiple} unstable modes, i.e., modes with amplification factors $|\sigma| > 1$. These comprise a small fraction of the overall spectrum, but nonetheless predict that there are modes of the linearized problem \eqref{eq:KdV-lin} that will grow over time. We limit our focus to the fastest growing of these unstable modes. For a given $(\alpha, \Delta t)$, these are the eigenvalues $\sigma$ with the largest absolute value, which correspond to the growth rates $\lambda$ with largest real part. If the linear theory is successful, these fastest modes will serve as 
reasonable estimates for the growth of the numerical instability in the corresponding simulations.

In Figure~\ref{fig:figure-2-evp-ivp}, we compare the behavior of the nonlinear numerical simulations with the exponential growth of the most unstable linear eigenmode. We find reasonable agreement in the unstable growth in the higher order SBDF schemes---SBDF3 and SBDF4---which clearly exhibit linear instabilities at early times. As for the remaining four schemes---SBDF1, SBDF2, RK222, and RK443---it is difficult to see any correspondence between the linear analysis and the unstable behavior in the nonlinear problem.

\begin{figure}[t]
    \centering
    \includegraphics{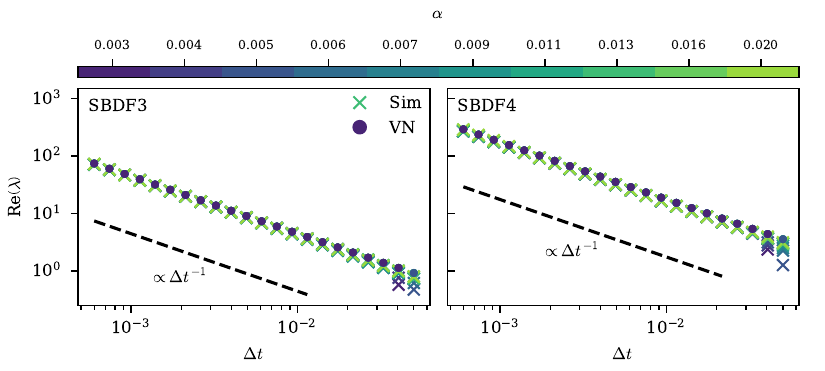}
    \caption{\textbf{Linear instability error growth rates measured in simulations (Sim) of solitons with SBDF3 and SBDF4 compared to rates predicted by von Neumann (VN) analysis}. Measured growth rates for both schemes exhibit good agreement with the rates predicted by taking the largest VN growth rates, with better agreement as $\Delta t \rightarrow 0$. Shown are (crosses) the measured exponential growth rates of the $L_2$ norm of the error in simulations with SBDF3 and SBDF4 and (circles) the real parts of the growth rates $\lambda$ inferred from the largest eigenvalues in the von Neumann eigenspectra. The former values were obtained from least squares fits of the $L_2$ data. Each point represents a measurement/computation for a particular ($\alpha, \Delta t$). The growth rate $\lambda$ is defined such that $L_2(u_{\rm sim} - u_{\rm ex}) \propto e^{{\rm Re}(\lambda) t}$. Both schemes exhibit ${\rm Re}(\lambda) \sim {\cal O}(\Delta t^{-1})$, indicating that reduction of the timestep size leads to faster unstable growth.}
    \label{fig:VN-sbdf3-4}
\end{figure}

Figure~\ref{fig:VN-sbdf3-4} shows that the linear analysis agrees well with SBDF3 and SBDF4 simulations across the full range of $\alpha$ and $\Delta t$. Intriguingly, while the growth rates in SBDF3 and SBDF4 have strong dependence on $\Delta t$, there is little sensitivity on the dispersion parameter $\alpha$ (in contrast to our observations of SBDF2 in Section~\ref{sec:sol}). The parameter study in Figure~\ref{fig:VN-sbdf3-4} also surprisingly finds that reduction in the timestep size $\Delta t$ corresponds to \textit{faster} growth rates, with ${\rm Re}(\lambda) \sim {\cal O}(\Delta t^{-1})$. This suggests that not only are SBDF3 and SBDF4 highly unstable for the solution of KdV, but their growth cannot be bounded via reduction of the timestep size. That is, SBDF3 and SBDF4 fail to satisfy the weak version of stability (which we call weakly unstable) where $|\sigma| \le 1 + {\cal O}(\Delta t)$ as $\Delta t \rightarrow 0$, which would allow one to ensure the solution remains \textit{bounded} up to a desired fixed time $T$ independent of timestep size, (e.g., \cite{Trefethen1994} or \cite{Deriaz2020}). 

We are more interested, however, in the failure of the linear analysis to account for the behavior of the four remaining schemes. The initial error levels in the first several iterations with SBDF1, SBDF2, RK222, and RK443 are consistently set by the temporal truncation error of each scheme. From then on, the error begins to exhibit growth much faster than predicted by the largest amplification factors found in the linear analysis. Inspection of snapshots of the numerical solution indicates that this rapid error growth corresponds to an increasing phase difference between the exact soliton and numerical solution. The apparent ``saturation'' of the error thus corresponds to the numerical and exact solution being completely out of phase with one another, except for whenever the numerical solution has migrated so far away within the periodic domain that it once again overlaps with the exact solution for a brief time. 

While a phase difference is the first noticeable indicator of this long-term process, other concurrent changes contribute to a gradual increase in the error level, including changes in the amplitude and shape of the soliton. We note that whereas SBDF1, SBDF2, and RK222 solitons increase in amplitude as they speed up, RK443 solitons decrease in amplitude and slow down, as solitons propagate at a speed proportional to their amplitude. Finally, after these developments have played out for some time, SBDF1, SBDF2, and RK222 experience a sharp blow-up, while RK443 avoids this fate.

\subsection{Accounting for nonlinear contributions to error growth}\label{sec:L-vs-NL}

For many of our timestepping schemes, the errors grow much more rapidly at early times than predicted by linear stability analysis, suggesting that nonlinear accumulation of error may be important very early in the simulations.
Therefore, we would like to understand and compare the linear and nonlinear contributions to the change in the norm of the perturbation, i.e., their contributions to $\partial_t \|v\|$. In the absence of nonlinearities, our von Neumann analysis identifies a spectrum of stable and unstable modes with corresponding rates $\lambda$ such that $\partial_t \|v\| = {\rm Re}(\lambda)\|v\|$. On the other hand, we neglected a nonlinear term of the form $-v\partial_x v$ in order to arrive at that result, which itself could induce a rate of change on the order of $\|v\partial_x v\|$. 

Since the error at early times is dominated by the truncation error of the timestepping scheme, an estimate for the size of the quadratic nonlinearity, at some fixed sufficiently early time $t_{\rm early}$, is $\|v\partial_x v\| \sim {\cal O}(\Delta t^{2s})$ as $\Delta t \rightarrow 0$, where $s$ is the order of accuracy. Note that while the first several steps with higher order SBDF schemes are dominated by the first step with SBDF1, we consider $t_{\rm early}$ to be taken when the truncation errors of the higher order SBDF scheme have taken precedence. We verify the validity of this estimate with direct computations of the size of this term throughout our ($\alpha$, $\Delta t$) parameter space in Appendix~\ref{sec:NL-arg}.

For the linear contribution, we use the upper bound ${\rm Re}(\lambda_{\rm max})\|v\|$, where $\lambda_{\rm max}$ is the eigenvalue with maximum real part. We compare ${\rm Re}(\lambda_{\rm max})\| v\|$ with $\|v\partial_x v\|$ in the limit as $\Delta t \rightarrow 0$ since this is the limit relevant for accuracy. The linear analysis typically holds when the timestep size is sufficiently small, but this may not be the case if the nonlinear term asymptotes to zero slower than the linear contributions as $\Delta t \rightarrow 0$. The earlier linearization \eqref{eq:KdV-lin} is only valid insofar as the linear term ${\rm Re}(\lambda_{\rm max})\| v\|$ dominates in this limit. Otherwise, if the balance is comparable or swings in the other direction, we must conclude that the linearization was a poor assumption. 

We obtain estimates for the size of ${\rm Re}(\lambda_{\rm max})$ from our full survey of the fastest growing von Neumann eigenmodes for the four schemes of interest (given in Appendix~\ref{sec:NL-arg}). These results along with those presented for SBDF3 and SBDF4 in Figure~\ref{fig:VN-sbdf3-4} allow us to compare linear and nonlinear contributions to error growth  in all schemes in Table~\ref{tbl:err-L-vs-NL}. The comparison finds that SBDF3 and SBDF4, where our linear analysis performs well, each pass the test of ${\rm Re}(\lambda_{\rm max})\|v\|$ dominating as $\Delta t \rightarrow 0$, whereas this is \textit{not} the case for the four schemes of interest, SBDF1, SBDF2, RK222, and RK443. This is a strong indication that linear analysis will not capture the behavior in these simulations, even at early times. Thus, we will henceforth work with the full nonlinear equation to understand the timestepping instability. 

\begin{table}%[h]
\caption{\textbf{Expected sizes of nonlinear error terms compared to growth rates from the VN analysis.} Here we summarize whether the linear theory is valid for each scheme by comparing the nonlinear and linear contributions to error growth at early times in the $\Delta t \rightarrow 0$ limit. The size of the nonlinear contribution is given by $\|v\partial_x v\|$ and an upper bound of the size of the linear contribution is given by ${\rm Re}(\lambda_{\rm max})\|v\|$, where $v = u_{\rm sim} - u_{\rm ex}$ is the error at a fixed early time and $\lambda_{\rm max}$ is the fastest VN growth rate. We take $\|\cdot\|$ to be the $L_2$ norm.}\label{tbl:err-L-vs-NL}
\begin{tabular*}{\tblwidth}{@{} L | L L L L @{} }
\toprule
 Scheme & $\|v \partial_x v\|$ & ${\rm Re}(\lambda_{\rm max})\|v\|$ & Comparison & Result \\
\midrule \midrule
 SBDF1 & ${\cal O}(\Delta t^2)$ & ${\cal O}(\Delta t^2)$ & Same order in $\Delta t$, neither dominates & Linear theory inadequate \\
 SBDF2 & ${\cal O}(\Delta t^4)$ & ${\cal O}(\Delta t^5)$ & Nonlinear term dominates as $\Delta t \rightarrow 0$ & Linear theory inadequate\\
 RK222 & ${\cal O}(\Delta t^4)$ & ${\cal O}(\Delta t^5)$ & Nonlinear term dominates as $\Delta t \rightarrow 0$ & Linear theory inadequate\\
 RK443 & ${\cal O}(\Delta t^6)$ & ${\cal O}(\Delta t^6)$ & Same order in $\Delta t$, neither dominates & Linear theory inadequate\\
 \midrule \midrule
 SBDF3 & ${\cal O}(\Delta t^6)$ & ${\cal O}(\Delta t^2)$ & Linear term dominates as $\Delta t \rightarrow 0$ & Linear theory successful\\
 SBDF4 & ${\cal O}(\Delta t^8)$ & ${\cal O}(\Delta t^3)$ & Linear term dominates as $\Delta t \rightarrow 0$ & Linear theory successful \\
\bottomrule
\end{tabular*}
\end{table}

\section{Multiple scales analysis}\label{sec:ms}

\subsection{Motivation of the multiple scales framework}

We turn our attention to the four schemes (SBDF1, SBDF2, RK222, and RK443) exhibiting behavior unexplained by the linear analysis in Section~\ref{sec:lin-stab}. While we have emphasized how the instability in question plays out over long time scales, it is just as critical to recognize that the behavior of the numerical solution on short time scales approximately tracks that of an exact soliton. Thus, there is a slow time scale on which the numerical solution experiences unstable error growth, and a fast time scale on which the solution's behavior is close to that of an analytic solution. 

Adopting this perspective, one might expect a \textit{multi-scale asymptotics} approach to be capable of yielding an approximation of the numerical solution that reproduces this separation of time scales. In this work, we will employ a multiple scales framework in the semi-discrete limit where the numerical solution is evolved by a timestepping scheme but remains continuous in space, as the spatial discretization does not impact the instability. 
An asymptotic framework has the further advantage of retaining the nonlinearities of the KdV equation, whereas the linear analysis does not explain the instability. 

Our approach is to define a small dimensionless parameter $\epsilon$ as the ratio of the timestep size to the propagation timescale across the soliton 
\begin{equation}\label{eq:ms-epsilon}
    \epsilon = \frac{\Delta t}{\tau_0} = \frac{\Delta t}{(x_{{\rm FWHM},0}/c_0)} = \frac{\Delta t \alpha^{-1/2} c_0^{3/2}}{4\ln{\left(1+\sqrt{2}\right)}}. %\propto \alpha^{-1/2} \Delta t.
\end{equation}
We take the propagation time scale $\tau_0$ to be the ratio of the initial full width at half maximum of the soliton to its initial propagation speed $c_0$. Thus, $\epsilon$ serves as a reference for the relative separation between discrete and continuum time scales in the problem.

We now aim to predict the nonlinear evolution of the numerical error in the limit of $\epsilon \rightarrow 0$ for a given timestepping scheme by allowing its solution to vary over multiple time scales depending on $\epsilon$. For an $s$-order accurate timestepping scheme, we introduce the $s+1$ time scales
\begin{equation}\label{eq:ms-times}
    t_0 = t, \quad t_1 = \epsilon t, \quad \dots, \quad t_s = \epsilon^s t, \quad t_{s+1} = \epsilon^{s+1}t.
\end{equation}
We then form a regular asymptotic expansion of the numerical solution $U^{(n)} \equiv U(x, t = n\Delta t)$ given by 
\begin{equation}\label{eq:ms-U}
    U^{(n)}(x) = \sum_{k=0}^{s+1}\epsilon^k U_k^{(n)}(x, t_{0}, \dots, t_{s+1}) + {\cal O}(\epsilon^{s+2}),
\end{equation}
where we have allowed each term in the expansion to possibly vary on each of the time scales introduced in \eqref{eq:ms-times}. 

So far we have used continuous variables to represent the multiple time scales, whereas the semi-discrete equations for a given timestepping scheme will consist of the solution at several time levels. To deal with this complication, we express the solution at each time level in terms of $U^{(n)}$ by performing Taylor expansions in $t$. This requires that we represent partial derivatives with respect to time $t$ in terms of the multiple time scales we have introduced, i.e., 
\begin{equation}\label{eq:ms-time-der-series}
    \frac{\partial}{\partial t} = %\frac{\partial t_0}{\partial t}\frac{\partial}{\partial t_0} + \frac{\partial t_1}{\partial t}\frac{\partial}{\partial t_1} + \dots =
    \sum_{m=0}^{s+1} \epsilon^m\frac{\partial}{\partial t_m} + {\cal O}(\epsilon^{s+2}).
\end{equation}
We therefore express the numerical solution at time level $t_{n+\nu} = t_n + \nu \Delta t$ wherever it appears in the semi-discrete timestepping equation as
\begin{equation}\label{eq:ms-taylor}
    U^{(n+\nu)}(x) = \sum_{k=0}^{s+1} \left[\epsilon^k U_k^{(n)}(x, t_{0}, \dots, t_{s+1})\right] + \nu\Delta t \left(\sum_{m=0}^{s+1}\left[\epsilon^m \frac{\partial}{\partial t_m}\right]\right)\left(\sum_{k=0}^{s+1} \left[\epsilon^k U_k^{(n)}(x, t_{0}, \dots, t_{s+1})\right]\right) + \dots.
\end{equation}

This approach allows us to seek an asymptotic approximation to $U^{(n)}(x)$ satisfying a system of continuous partial differential equations %, not in $(x,t)$, but rather 
in $(x,t_0, t_1, \dots)$. An alternative formulation---in light of the method introduced in \cite{Torng1960} and developed further in later works (e.g., \citealt{Subramanian1979}, \cite{vanHorssen2009})---would be to introduce multiple \textit{discrete} iterative time scales, e.g., $n_0 = n, n_1 = \epsilon n, \dots$, which we would expect to lead to a corresponding system of partial \textit{difference} equations. Instead, our approach most closely resembles that of \cite{Hoppensteadt1977} which similarly bridges from a discrete system to a continuous one by way of Taylor expansions. Since our instability is related to timestepping errors, we found the latter approach more amenable as it  mirrors standard methods for computing the temporal truncation error. This construction of our multiple scales framework will serve to approximate the long-term accumulation of timestepping errors over many iterations.

\subsection{Fast-time dynamics of the numerical solution}

We begin by matching terms in the semi-discrete timestepping equations of order ${\cal O}(\epsilon^0)$. In multi-step schemes the matching procedure is made possible by expressing terms in the form described in~\eqref{eq:ms-taylor}, but multi-stage schemes have the additional complication of having intermediate stages to contend with. It will be sufficient to %facilitate the analysis by assuming 
assume that each intermediate stage has its own arbitrary regular expansion in $\epsilon$. %, which we won't need to determine.
With these considerations, we will suppress the superscript in the series expansion of $U^{(n)}$ in \eqref{eq:ms-U}, so that we have $U_0 \equiv U_0^{(n)}(x), U_1 \equiv U_1^{(n)}(x),$ and so on.

%We begin by matching terms of order ${\cal O}(\epsilon^0)$, after all necessary substitutions have been implemented. This corresponds to terms both at leading order in the numerical solution $U$ and in time $t$. 
In each of the four timestepping schemes of interest, we find the leading order problem
\begin{equation}\label{eq:ms-leading-order}
    \begin{aligned}
    {\cal O}(\epsilon^0) \quad &\partial_{t_0} U_0 + \alpha \partial_{x}^3 U_0 + U_0 \partial_x U_0 = 0, \quad -\frac{L}{2} \le x \le \frac{L}{2}, \quad t \ge 0,\\
    &U_0|_{x=-L/2} = U_0|_{x=L/2},\quad \partial_x U_0|_{x=-L/2} = \partial_x U_0|_{x=L/2},\quad \partial_x^2 U_0|_{x=-L/2} = \partial_x^2 U_0|_{x=L/2}.
    \end{aligned}
\end{equation}
This states that the leading order numerical solution $U_0$ satisfies the KdV equation with respect to the fast time $t_0$. We suppose in~\eqref{eq:ms-leading-order} that $U_0$ is periodic in $x$ to be consistent with the periodic boundary conditions in our numerical simulations (and later we will assume higher order terms in the expansion of $U^{(n)}$ are likewise periodic). 

Analytic solutions to the periodic KdV equation are wave train solutions called cnoidal waves, which are expressed in terms of the Jacobi elliptic functions ${\rm cn}(x; m)$ with parameter $0 < m < 1$ setting the period. However, these solutions reduce to the soliton~\eqref{eq:soliton} when the period is much longer than the soliton width, since ${\rm cn}(x; m) \rightarrow {\rm sech}(x)$ as $m \rightarrow 1$ (Chapter 16 of \citealt{Abramowitz1964}). We will therefore consider solutions to~\eqref{eq:ms-leading-order} to be in the limiting case of $m = 1$ as this most closely resembles the length scales in our simulations. 

An additional consideration for the expression of $U_0$ is that~\eqref{eq:ms-leading-order} has a Galilean invariance, 
so the integral of its solutions $U_0$ is conserved over time. Likewise, each of our multi-step and multi-stage temporal discretizations of the KdV equation will conserve this density, as this is guaranteed by the order conditions of each scheme. We therefore arrive at the more general expression for $U_0$ satisfying~\eqref{eq:ms-leading-order} given by
\begin{equation}\label{eq:ms-U0-general}
    U_0(x, t_0, t_1, \dots, t_{s+1}) = u_{\infty}(c) + 3c\,{\rm sech}{\left(\sqrt{\frac{c}{4\alpha}}(x-(c + u_{\infty}(c))t_0\right)^2},
\end{equation}
where we have leveraged the conserved quantity $\int_{-L/2}^{L/2} U_0\, dx$ to define the shift parameter $u_{\infty}$ given by 
\begin{equation}\label{eq:ms-uinf}
    u_{\infty} \equiv u_{\infty}(c) =\frac{12\alpha^{1/2}}{L}\left(c_0^{1/2}{\rm tanh}{\left(\frac{L}{2}\sqrt{\frac{c_0}{4\alpha}}\right)} - c^{1/2}{\rm tanh}{\left(\frac{L}{2}\sqrt{\frac{c}{4\alpha}}\right)}\right). %\approx \frac{12\alpha^{1/2}}{L}\left(c_0^{1/2} - c^{1/2}\right).
\end{equation}
Note that~\eqref{eq:ms-uinf} distinguishes between a general value of $c$ and its initial value $c_0$ established by the initial condition. The significance of this may not be apparent until we express that
\begin{equation}\label{eq:ms-c}
    c \equiv c(t_1, \dots, t_s, t_{s+1}),
\end{equation}
emphasizing that $c$ in~\eqref{eq:ms-U0-general} and \eqref{eq:ms-uinf} must be fixed with respect to $x$ and $t_0$, but has not yet been precluded from varying \textit{with respect to the slow time scales}. This observation will prove important in the work that follows.

\subsection{Slow-time evolution of the numerical solution}\label{sec:ms-slow}

For $s$-order accurate SBDF or ERK+DIRK schemes, we encounter consistent patterns when matching terms of order ${\cal O}\left(\epsilon^m\right)$ where $m < s$. At these orders we find equations of the form
\begin{equation}\label{eq:ms-interm-order}
    {\cal O}(\epsilon^{m}) \quad \partial_{t_0} U_m + \alpha \partial_{x}^3 U_m + \partial_x(U_0 U_m) = -\partial_{t_m} U_0 + g_m(U_1, \dots, U_{m-1}), \quad 0 < m < s
\end{equation}
where $g_m$ is in general a nonlinear function of $U_1, \dots, U_{m-1}$ and their derivatives---however $g_m$ does not depend on $U_0$. The right hand side of \eqref{eq:ms-interm-order} also includes %importantly consists of 
the term $-\partial_{t_m} U_0$ capturing the slow time dependence of $U_0$ on $t_m =  \epsilon^m t$. Under the change of variables $\xi = x - ct_0$, the left hand side is the adjoint of the KdV operator acting on $U_m$. %The left hand side on the other hand resembles the KdV equation (in $(x, t_0)$) written for $U_1$ with respect to the fast time scale $t_0$, except the advection term is linear in $U_1$ and contains $U_0$. We recognize that the left hand side is in fact formally adjoint to the KdV equation when both are considered under the change of variables $\xi = x - ct_0$. 

Since the problem will be periodic with respect to $\xi$, we can leverage the adjoint structure by applying this change of variables to~\eqref{eq:ms-interm-order}, multiplying through by $U_0$, and integrating by parts in $\xi$. After applying periodic boundary conditions all left hand side terms vanish, leaving us with a \textit{solvability condition} of the form
\begin{equation}\label{eq:ms-sc-m}
    %\int_{-L/2}^{L/2} U_0\partial_{t_m} U_0\, d\xi  %= 
    \frac{1}{2} \partial_{t_m}{\left( \int_{-L/2}^{L/2} U_0^2\,d\xi \right)}
    = \int_{-L/2}^{L/2} U_0 g_m(U_1, \dots, U_{m-1})d\xi.
\end{equation}
Since $g_m$ does not depend on $U_0$, we can satisfy this equation %are able to simplify this equation 
by taking $U_k = 0$ for $0 < k < m$. Then $g_m$ and thereby the right hand side of~\eqref{eq:ms-sc-m} can be taken to be zero, and consequently we are able to conclude that
%\begin{equation}\label{eq:trivial-sc}
    %\frac{1}{2} \partial_{t_m} \left( \int_{-L/2}^{L/2} U_0^2\,d\xi \right)= 0,
%\end{equation}
%which requires that 
$U_0$ has \textit{no dependence} on this order's slow time scale, $t_m$. 

However, a small but significant change occurs when we match terms at ${\cal O}(\epsilon^s)$ and beyond. When the order in $\epsilon$ matches the order of accuracy of the timestepping scheme, we now find equations of the form
\begin{equation}\label{eq:ms-order-s}
    {\cal O}(\epsilon^{s}) \quad \partial_{t_0} U_s + \alpha \partial_{x}^3 U_s + \partial_x(U_0 U_s) = -\partial_{t_s} U_0 + g_s(U_0, U_1, \dots, U_{s-1}).
\end{equation}
The right hand side of~\eqref{eq:ms-order-s} differs crucially from \eqref{eq:ms-interm-order} in that $g_s$ now consists of terms depending on $U_0$ and its derivatives---which cannot be neglected. The solvability equation then takes the form %We are consequently forced to accept a solvability condition of the form
\begin{equation}\label{eq:ms-sc-s}
    \frac{1}{2}\partial_{t_s}{\left(\int_{-L/2}^{L/2}U_0^2\, d\xi\right)} = \int_{-L/2}^{L/2}U_0 g_s(U_0)\,d\xi.
\end{equation}
Here we emphasize that $g_s(U_0)$ only depends on the leading order term $U_0$, as $U_1, \dots, U_{s-1}$ are set to zero or left out of the original expansion, and all $U_s$ dependence has been removed via integration by parts. We therefore expect~\eqref{eq:ms-sc-s} to give us a \textit{non-trivial} solvability condition which will provide a description of how $U_0$ depends on the slow time scale $t_s$. %As we have hinted at prior to now, the 
This non-trivial solvability condition will tell us how the parameter $c$ varies with respect to $t_s$. In fact, the condition \eqref{eq:ms-sc-s} typically \textit{reduces to a first-order differential equation} for $c$ in the slow time scale $t_s$.

For some of our timestepping schemes, we find that the right hand side of~\eqref{eq:ms-sc-s} vanishes due to parity considerations. In particular, SBDF2 and RK222---the two schemes we consider that are second-order accurate---admit $g_s$ that are odd functions with an anti-symmetry about $\xi = 0$. This is a consequence of the alternating parity of $U_0$ (which is even in $\xi$) and its derivatives: when $s$ is an even positive integer, $g_s(U_0)$ consists of terms with only an odd number of derivatives of $U_0$, whereas when $s$ is odd the opposite occurs. Therefore, when $s=2$ as in SBDF2 and RK222, we again obtain a trivial solvability condition at ${\cal O}(\epsilon^s)$ and must go one order further to ${\cal O}(\epsilon^{s+1})$ to obtain a nontrivial condition, and thereby a first-order differential equation for $c$ with respect to the slow time scale $t_{s+1}$. Note, however, that while we find that $\partial_{t_s} U_0 = 0$, we can not set $U_s = 0$ as this would impose inconsistent constraints on $U_0$. Fortunately, this does not prevent us from obtaining the condition at next order. 

\section{Results}\label{sec:results}

\subsection{Predictions of the multiple scales analysis}\label{sec:ms-solvability-conditions}

\subsubsection{SBDF1 (Forward-Backward Euler)}\label{sec:results-sbdf1}
When evolving KdV solitons with the first-order accurate IMEX scheme SBDF1, the multiple scales analysis finds the expected fast-time behavior given in~\eqref{eq:ms-leading-order}, and encounters a nontrivial result at order ${\cal O}(\epsilon^1)$. Matching terms at this order gives the equation
\begin{equation}\label{eq:sbdf1-e1}
    {\cal O}(\epsilon^1) \quad \partial_{t_0} U_1 + \alpha\partial_x^3 U_1 + \partial_x(U_0 U_1) = - \partial_{t_1} U_0 + \tau_0\partial_{t_0}{\left(-\frac{\alpha}{2}\partial_x^3U_0 + \frac{1}{2}U_0\partial_xU_0\right)},
\end{equation}
which admits a nontrivial solvability condition of the form
\begin{equation}\label{eq:sbdf1-sc}
    \frac{1}{2}\partial_{t_1}{\left(\int_{-L/2}^{L/2} U_0^2\,d\xi\right)} = \int_{-L/2}^{L/2} U_0\tau_0\partial_{t_0}{\left(-\frac{\alpha}{2}\partial_x^3U_0 + \frac{1}{2}U_0\partial_xU_0\right)}\, d\xi.
\end{equation}
where $\xi = x - ct_0$, such that $\partial_x = \partial_{\xi}$ and $\partial_{t_0} = -c\partial_{\xi}$. When we assume that $U_0$ has the form of the soliton expressed in~\eqref{eq:ms-U0-general}, we find that~\eqref{eq:sbdf1-sc} reduces to an ordinary differential equation (ODE) for the parameter $c$ of the form
\begin{equation}\label{eq:sbdf1-c-ode}
    \frac{\partial c}{\partial t_1} = \frac{{\rm RHS}_{\rm SBDF1}(c; \alpha, L, c_0)}{{\rm LHS}(c;\alpha, L, c_0)}.
\end{equation}
When working on the finite periodic domain, the terms in~\eqref{eq:sbdf1-c-ode} consist of nonlinear hyperbolic trigonometric functions of $c$ which must be integrated numerically. Expressions of the numerator and denominator in~\eqref{eq:sbdf1-c-ode} are provided in the Supplemental Material.

However, we can obtain a corresponding condition that can be solved {\it analytically} if we consider the limit of an infinite domain by taking $L \rightarrow \infty$ (such that $u_{\infty} \rightarrow 0$) and instead impose decay boundary conditions on $U_0$ and its first and second spatial derivatives. In this case, the solvability condition \eqref{eq:sbdf1-c-ode} reduces to the simple separable ODE
\begin{equation}\label{eq:sbdf1-ode-related}
    \frac{\partial c}{\partial t_1} = \frac{34}{105}\frac{\tau_0}{\alpha}c^4.
\end{equation}
When equipped with the initial condition $c(0) = c_0$, \eqref{eq:sbdf1-ode-related} admits the \textit{power law} solution
\begin{equation}\label{eq:sbdf1-PL-related}
    c(t_1) = \left(c_0^{-3} - \frac{34}{35}\frac{\tau_0}{\alpha}t_1\right)^{-1/3} \quad \Leftrightarrow \quad c(t) = \left(c_0^{-3} - \frac{34}{35}\frac{\Delta t}{\alpha}t\right)^{-1/3},
\end{equation}
where we used $t_1 = \epsilon t = (\Delta t/\tau_0) t$ to express $c$ as a function of $t$. This predicts  that the parameter $c$---and hence the amplitude of the soliton---will exhibit unstable growth via a nonlinear power law in time. Hence, the numerical solution is predicted to diverge to $+\infty$ at the finite time $t = (34/35)c_0^{-3}(\alpha/\Delta t)$. 

We expect that the result in~\eqref{eq:sbdf1-PL-related} will be a good approximation to the solution of~\eqref{eq:sbdf1-c-ode} in the limit that the soliton length scale set by $\alpha$ is much smaller than the domain scale set by $L$. It is not feasible to arrive at a similar closed-form prediction when the problem is posed on the finite periodic interval, which more accurately captures our numerical setup. We therefore numerically integrate the condition \eqref{eq:sbdf1-c-ode} to compute numerical predictions, while using \eqref{eq:sbdf1-PL-related} to facilitate interpretations.

\subsubsection{SBDF2 and RK222}\label{sec:results-sbdf2-rk222}

By matching terms at ${\cal O}(\epsilon^0)$ and ${\cal O}(\epsilon^1)$, we find that the second-order accurate IMEX schemes SBDF2 and RK222 have the expected fast-time behavior described by~\eqref{eq:ms-leading-order} and an equation of the form~\eqref{eq:ms-interm-order} admitting the trivial solvability condition that $U_0$ has no $t_1$ dependence. As described in Section~\ref{sec:ms-slow}, at order ${\cal O}(\epsilon^2)$, both schemes have equations of the form~\eqref{eq:ms-order-s} where extra terms depending on $U_0$ remain. However, this solvability condition is trivially satisfied due to the parity of $U_0$ as a function of $\xi$. Therefore, our work up to ${\cal O}(\epsilon^2)$ gives that $U_1 = 0, \partial_{t_1} U_0 = 0$, and $\partial_{t_2} U_0$ = 0. 

Taking SBDF2 for example, we then find the ${\cal O}(\epsilon^3)$ equation
\begin{equation}\label{eq:sbdf2-e3}
    \partial_{t_0}U_3 + \alpha\partial_x^3 U_3 + \partial_x(U_0 U_3) = -\partial_{t_3}U_0 + \tau_0^3\partial_{t_0}^3{\left(\frac{\alpha}{4}\partial_x^3U_0 - \frac{3}{4}U_0\partial_xU_0\right)}.
\end{equation}
This corresponds to the nontrivial solvability condition for SBDF2
\begin{equation}\label{eq:sbdf2-e3-sc}
    \frac{1}{2}\partial_{t_3}{\left(\int_{-L/2}^{L/2}U_0^2\,d\xi\right)} = \int_{-L/2}^{L/2} U_0\tau_0^3\partial_{t_0}^3{\left(\frac{\alpha}{4}\partial_x^3U_0 - \frac{3}{4}U_0\partial_xU_0\right)}\, d\xi,
\end{equation}
where the right hand side integrand consists of terms with an even number of derivatives of $U_0$, such that the right hand side does not vanish. For RK222, the equations at ${\cal O}(\epsilon^3)$ and solvability condition are analogous but significantly more complicated, so we provide them in the Supplemental Material. 

As we did for SBDF1 in Section~\ref{sec:results-sbdf1}, when we solve the problem over a finite periodic interval, we obtain ODEs for $c$ that can be integrated numerically (Supplemental Material)---but this time the independent variable is the slow time scale $t_3$. Alternatively, for an infinite domain ($L \rightarrow \infty$), the respective nontrivial solvability conditions are
\begin{equation}\label{eq:sbdf2-ode-related}
    \frac{\partial c}{\partial t_3} = \frac{43}{105}\frac{\tau_0^3}{\alpha^2}c^7 \quad \text{(SBDF2)}
\end{equation}
and
\begin{equation}\label{eq:rk222-ode-related}
    \frac{\partial c}{\partial t_3} = \frac{355045-245436\sqrt{2}}{30030}\frac{\tau_0^3}{\alpha^2}c^7 \quad \text{(RK222).}
\end{equation}
With $c(0) = c_0$, these ODEs again correspond to power law solutions, respectively given for SBDF2 and RK222 by
\begin{equation}\label{eq:sbdf2-PL-related}
    c(t_3) = \left(c_0^{-6} -\frac{86}{35}\frac{\tau_0^3}{\alpha^2}t_3\right)^{-1/6} \quad \Leftrightarrow \quad c(t) = \left(c_0^{-6} -\frac{86}{35}\frac{\Delta t^3}{\alpha^2}t\right)^{-1/6}
\end{equation}
and
\begin{equation}\label{eq:rk222-PL-related}
    c(t_3) = \left(c_0^{-6} -\frac{355045-245436\sqrt{2}}{5005}\frac{\tau_0^3}{\alpha^2}t_3\right)^{-1/6} \Leftrightarrow c(t) = \left(c_0^{-6} -\frac{355045-245436\sqrt{2}}{5005}\frac{\Delta t^3}{\alpha^2}t\right)^{-1/6}.
\end{equation}
We note that this $L \rightarrow \infty$ approximation predicts these schemes to have blow-up times scaling like $\alpha^2 \Delta t^{-3}$. This again captures the unstable growth in the soliton, but according to a power law that blows up on slower time scales and with a sharper singularity. 

\subsubsection{RK443}\label{sec:results-rk443}

As expected, we obtain a nontrivial solvability condition for the third-order accurate IMEX scheme RK443 when matching terms at ${\cal O}(\epsilon^3)$. The equation is quite lengthy and is provided in the Supplemental Material. Once again, when $U_0$ is posed as~\eqref{eq:ms-U0-general} over the finite periodic interval, this results in an ODE for $c$ with respect to the slow time scale $t_3$ that can be integrated numerically.

In the limit of an infinite domain ($L \rightarrow \infty$), the nontrivial solvability condition reduces to
\begin{equation}\label{eq:rk443-ode-related}
    \frac{\partial c}{\partial t_3} = -\frac{77069}{180180}\frac{\tau_0^3}{\alpha^2}c^7.
\end{equation}
Unlike the prior timestepping schemes, the terms in \eqref{eq:rk443-ode-related} are opposite in sign. The corresponding power law solution with $c(0) = c_0$ is given by
\begin{equation}\label{eq:rk443-PL-related}
    c(t_3) = \left(c_0^{-6} + \frac{77069}{30030}\frac{\tau_0^3}{\alpha^2}t_3\right)^{-1/6} \quad \Leftrightarrow \quad c(t) = \left(c_0^{-6} + \frac{77069}{30030}\frac{\Delta t^3}{\alpha^2}t\right)^{-1/6}.
\end{equation}
The consequence of the difference in sign in \eqref{eq:rk443-ode-related} is that the power law now \textit{decays} as $t\rightarrow \infty$, rather than blowing up. Consequently, RK443 differs from SBDF1, SBDF2, and RK222 in that, rather than growing in amplitude, decreasing in width, and gaining speed, its solitons will decay in amplitude, grow wider, and slow down. While in some sense opposite behaviors, we would still expect the norm of the numerical error to appear similar across all schemes, as we indeed saw earlier in Figure~\ref{fig:figure-2-evp-ivp}. The reason for this similarity is because each scheme is predicted to experience a phase shift from the exact solution due to the slow-time evolution of $c$. The schemes differ in that RK443 does not experience late-time blow-up. Rather, the multiple scales analysis predicts that RK443 is stable insofar as its error growth with respect to the exact solution is bounded, as the norm of the numerical solution is predicted to decay to zero as $t \rightarrow \infty$.

\subsubsection{Summary of the solvability conditions}

\begin{table}[b]
\caption{\textbf{Predictions from nontrivial solvability conditions in the multiple scales analysis}. The multiple scales analysis predicts that the respective numerical solutions of each scheme will evolve on a slow time scale related to an odd order in $\epsilon$. RK443 differs from SBDF1, SBDF2, and RK222 in that its solution is expected to exhibit decay, rather than blow-up, but all four schemes yield solutions that follow similar nonlinear power laws in time. The decay time for RK443 is chosen to be the time such that $L_2(U(x,T_{\rm decay})) = 0.9 L_2(U(x,0))$. To obtain $T_{{\rm decay}}$ for another decay fraction $0 < \phi < 1$, one would replace 0.9 by $\varphi$.
}
\label{tbl:ms-results}
\begin{tabular*}{\tblwidth}{@{} L | L L L @{} }
\toprule
 Scheme & Order of accuracy in $\Delta t$ & Order of slow-time dependence in $\epsilon$ & End behavior for solitons posed over $x \in (-\infty,\infty)$ \\ 
\midrule \midrule
SBDF1 & $\Delta t$ & $\epsilon$ & $T_{\rm blowup} = \frac{35}{34}c_0^{-3}\frac{\alpha}{\Delta t}$\\
SBDF2 & $\Delta t^2$ & $\epsilon^3$ & $T_{\rm blowup} = \frac{35}{86}c_0^{-6}\frac{\alpha^2}{\Delta t^3}$\\
RK222 & $\Delta t^2$ & $\epsilon^3$ & $T_{\rm blowup} = \frac{5005}{355045 - 245436\sqrt{2}}c_0^{-6}\frac{\alpha^2}{\Delta t^3}$ \\
RK443 & $\Delta t^3$ & $\epsilon^3$ & $T_{\rm decay} = \left(0.9^{-4} - 1\right)\frac{30030}{77069}c_0^{-6}\frac{\alpha^2}{\Delta t^3}$\\
\bottomrule
\end{tabular*}
\end{table}

Table~\ref{tbl:ms-results} summarizes the orders in $\epsilon$ at which each of these ODEs were recovered. Also, while the behavior is slightly more complicated when we consider a finite periodic interval, Table~\ref{tbl:ms-results} expresses the expected behavior of the solution as predicted in the more tractable $L \rightarrow \infty$ limit. These provide fairly accurate estimates for the behavior observed in our simulations and facilitate quick comparisons between each scheme.

In each of the four IMEX schemes of interest, we applied the multiple scales analysis up to a sufficiently high order in $\epsilon$ to obtain a nontrivial solvability condition. These conditions each resulted in an ODE which can be solved to find how the soliton parameter $c$, and thereby $U_0$, varies with respect to one of the slow time scales. When paired with the description of how $U_0$ behaves on fast time scales given by~\eqref{eq:ms-leading-order}, one then obtains a leading-order approximation to the numerical solution as it varies in $x$, the fast time $t_0$, and the identified slow time variable.

\subsection{Numerical results}\label{sec:ms-results}
As described in Section~\ref{sec:intro}, in this work we simulated KdV solitons of the initial form~\eqref{eq:soliton} over the finite periodic interval $-L/2 \le x \le L/2$ for logarithmically sampled dispersion parameters $0.003 \le \alpha \le 0.02$ and fixed timestep sizes $0.0006 \le \Delta t \le 0.05$. In this section we report on simulations performed under a real Fourier basis with $N = 512$ basis functions, though we emphasize that the numerical instability is insensitive to the choice of basis and resolution. 

\subsubsection{Blow-up times}\label{sec:blowups}
\begin{figure}[b]
    \centering
    \includegraphics[width = \textwidth]{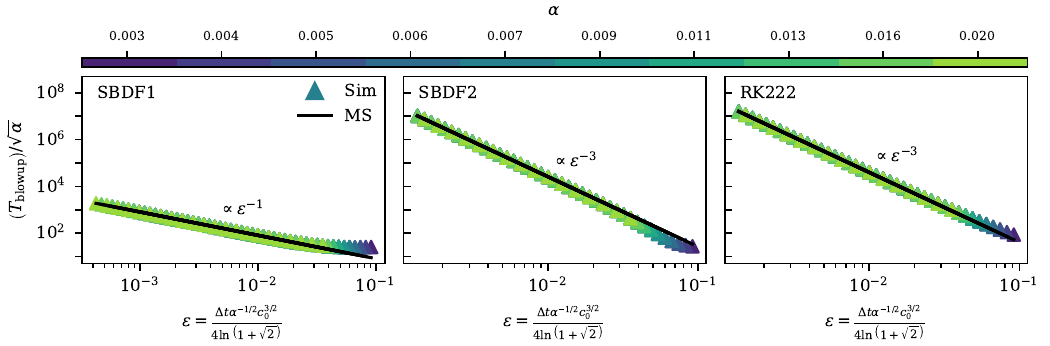}
    \caption{\textbf{Measured and predicted soliton blow-up times $T_{\rm blowup}$ in SBDF1, SBDF2, and RK222 simulations}. The accuracy of the blow-up time prediction improves as $\epsilon \rightarrow 0$. Blow-up times $T_{\rm blowup}$ in SBDF1, SBDF2, and RK222 are plotted against $\epsilon \propto \alpha^{-1/2} \Delta t$ as defined in \eqref{eq:ms-epsilon}. We plot $T_{\rm blowup}$ values measured in simulations (Sim) and those predicted by the multiple scales analysis (MS) assuming the solitons are evolved over a finite periodic domain. }
    \label{fig:t-ends-group}
\end{figure}

The simulations of KdV solitons ran using SBDF1, SBDF2, and RK222 exhibit finite-time blow-up. Figure~\ref{fig:t-ends-group} plots the blow-up times measured in our simulations against those predicted by the multiple scales analysis assuming the solitons are evolved over a finite periodic domain. As $T_{\rm blowup}$ scales like $\alpha^{1/2}$ multiplied by a power of $\epsilon$ in the $L \rightarrow \infty$ limit, we plot $T_{\rm blowup}\alpha^{-1/2}$ as a function of $\epsilon$. Although our solitons are evolved over a finite domain with $L = 10$, we expect this to work fairly well since we choose $\alpha$ such that the soliton widths are much smaller than $L$ (as discussed in Section~\ref{sec:sol}). 

Figure~\ref{fig:t-ends-group} finds good agreement between the measured and predicted $T_{\rm blowup}$. The magnitude and scaling with $\epsilon$ match our analytic calculations, particularly in the $\epsilon \rightarrow 0$ limit. We also note that at fixed $\epsilon$, we detect little to no variation in $T_{\rm blowup}\alpha^{-1/2}$ with respect to $\alpha$, verifying the $\alpha^{-1/2}$ normalization factor. 

\begin{figure}[h]
    \centering
    \includegraphics[width = \textwidth]{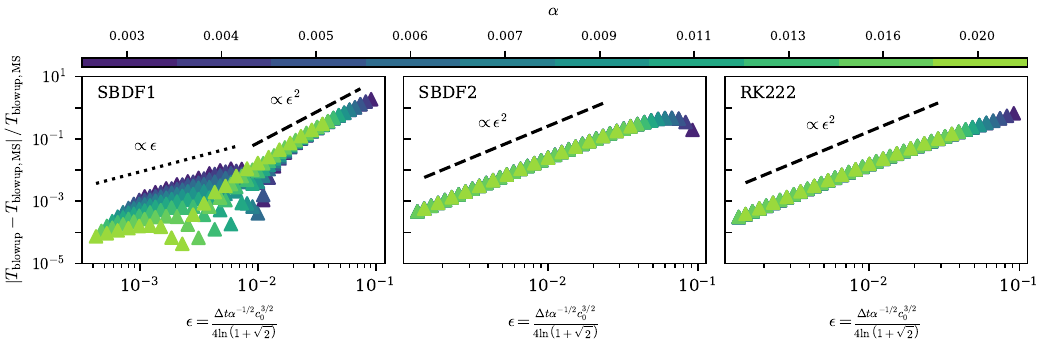}
    \caption{\textbf{Relative errors between measured blow-up times $T_{\rm blowup}$ and predicted blow-up times $T_{\rm blowup, MS}$ in SBDF1, SBDF2, and RK222 simulations}. 
    The relative error in the blow-up time prediction converges in the $\epsilon \rightarrow 0$ limit. As $\epsilon \rightarrow 0$, the errors decay as $\epsilon^2$, except for the SBDF1 timestepping scheme. In that case, we find errors instead decrease like $\epsilon$ for $\epsilon \lesssim 10^{-2}$.
    }
    \label{fig:t-ends-rel-err-group}
\end{figure}

In Figure~\ref{fig:t-ends-rel-err-group} we plot the relative errors in our predictions. %We continue to inspect the precision of the multiple scales analysis in Figure~\ref{fig:t-ends-rel-err-group} where we compute and plot the relative errors in our predictions,
In each scheme, the analytic predictions improve as $\epsilon \rightarrow 0$, with the relative errors typically converging $\propto \epsilon^2$. This is likely because $T_{\rm blowup, MS}$, the multiple scales analysis prediction of the blow-up time, is determined to some leading-order power of $\epsilon$ and lacks a correction at next-highest order due to parity considerations, i.e., the first nonzero corrections come two orders higher in $\epsilon$. 

The behavior of SBDF1 is an exception, as the convergence of its relative errors slows down to $\propto \epsilon$ at smaller values of $\epsilon$ with some dependence on $\alpha$. The switch in the convergence rate coincides with a change in sign of the relative error, suggesting the leading-order multiple scales prediction is superseded by another correction of opposite sign and with a less favorable convergence rate. Relatedly, our simulations of solitons with SBDF1 appear to develop slight asymmetries, mainly in the form of a spurious negative tail in the wake of the soliton, which can manifest over the long time scales of the simulation. This feature could violate the parity argument associated with the $\propto \epsilon^2$ convergence, suggesting the correction to $T_{\rm blowup, MS}$ would instead come at next order in $\epsilon$. %, especially in the limit of small $\epsilon$ where such asymmetries may have more time to accumulate. 
We have also noticed that a similar asymmetric tail sometimes develops in our higher-order simulations as well, but only in the final few iterations. Perhaps because this feature is shorter-lived than in SBDF1 simulations, this does not appear to have any appreciable effect on the convergence of the predictions in the higher-order schemes, at least for the range of $\epsilon$ considered in this work. The source of these asymmetries is an open question, as they are not explained by our leading-order predictions from the multiple scales analysis. Nevertheless, in all cases, we find the multiple scales prediction converges to the blow-up time of the simulations as $\epsilon \rightarrow 0$. 

\subsubsection{Decay time scales in RK443}\label{sec:decays}

\begin{figure}[h]
    \centering
    \includegraphics[width = \textwidth]{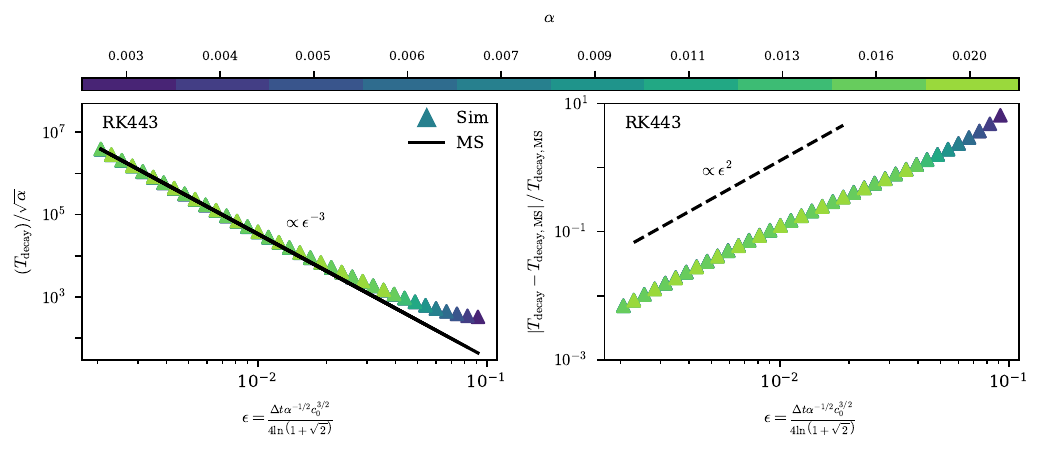}
    \caption{\textbf{Measured decay times $T_{\rm decay}$ and predicted decay times $T_{\rm decay, MS}$ in RK443 simulations and corresponding relative errors}. The decay time prediction for RK443 converges in the $\epsilon \rightarrow 0$ limit. (Left) Decay times in RK443 are plotted against $\epsilon \propto \alpha^{-1/2}\Delta t$ as defined in \eqref{eq:ms-epsilon}. We again plot values measured in simulations (Sim) and those predicted by the multiple scales analysis (MS) assuming the solitons are evolved over a finite periodic domain. 
    (Right) Corresponding absolute relative errors are plotted with respect to $\epsilon$, which we find are $\propto \epsilon^2$ in the $\epsilon \rightarrow 0$ limit.}
    \label{fig:t-ends-rk443}
\end{figure}

Whereas solitons in SBDF1, SBDF2, and RK222 are unstable, those in RK443 experience a gradual spurious decay, increasing in width while losing speed and amplitude such that the error is bounded as $t \rightarrow \infty$. Thus, rather than blow-up time, we calculate the time for the norm of the soliton to decrease by 10\%. Figure~\ref{fig:t-ends-rk443} plots the decay times $T_{\rm decay}$ such that 
$L_2(U(x,T_{\rm decay})) = 0.9 L_2(U(x,0))$. Figure~\ref{fig:t-ends-rk443} also includes the multiple scales prediction of this decay time scale, and the relative error in the prediction. 

We again find that the multiple scales analysis is effective in predicting the relevant slow time scale on which the numerical solution evolves. The successful analysis of RK443 is an important milestone in this work. This example demonstrates how the multiple scales framework is able to distinguish between cases when the solution will remain bounded, versus those in which the nonlinear evolution of the temporal error will be unstable and result in finite-time blow-up.

Moreover, although the numerical solution in RK443 does deviate from the exact analytic solution, the analysis predicts the solution to be stable for sufficiently small $\epsilon$. If one is choosing between these four IMEX methods to evolve solitons on arbitrarily long time scales, the analysis shows that RK443 is preferable for its prolonged stability. On the other hand, RK443 shares in common with SBDF2 and RK222 the $\alpha^2 \Delta t^{-3}$ scaling of their respective decay/blow-up time scales. Hence, if one only needs to propagate solitons for a set amount of time, then the lower order methods of SBDF2 and RK222 may offer more efficient computation for a given acceptable level of error. 

\subsubsection{Evolution of the error at intermediate times}

\begin{figure}[b]
    \centering
    \includegraphics[width = \textwidth]{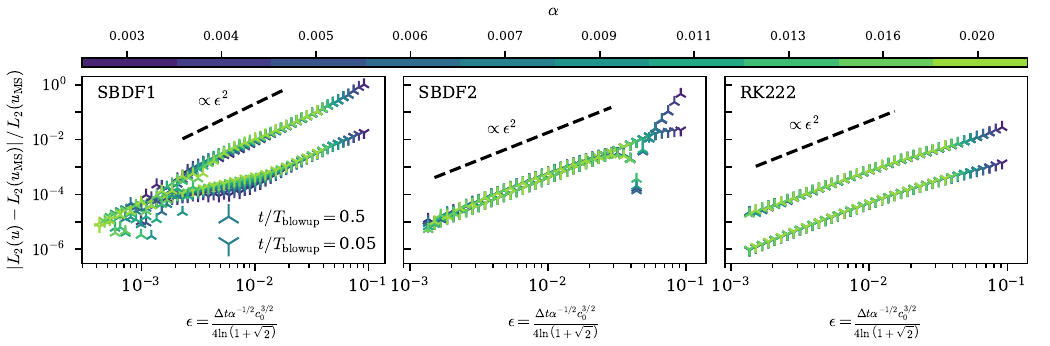}
    \caption{\textbf{Relative errors between measured norms $L_2(u)$ and predicted norms $L_2(u_{\rm MS})$ at intermediate times in SBDF1, SBDF2, and RK222 simulations}. Predictions of the $L_2$ norm of the numerical solution converge in the $\epsilon \rightarrow 0$ limit. Errors are considered at two intermediate times of $t = 0.05 T_{\rm blowup}$ and $t = 0.5 T_{\rm blowup}$, where $T_{\rm blowup}$ is either the simulation blow-up time or the predicted blow-up time for $L_2(u)$ and $L_2(u_{\rm MS})$, respectively.
    We typically find that errors are higher at later times, though this is difficult to see with SBDF2, possibly because the analysis neglects the error introduced by taking the first step with SBDF1. As $\epsilon \rightarrow 0$, the errors generally decrease as
    $\epsilon^2$. Errors for SBDF1 at small $\epsilon$ are an exception, which we again attribute to asymmetries not captured by the leading-order multiple scales prediction.}
    \label{fig:err-wrt-group}
\end{figure}

\begin{figure}[t]
    \centering
    \includegraphics[width = 0.5\textwidth]{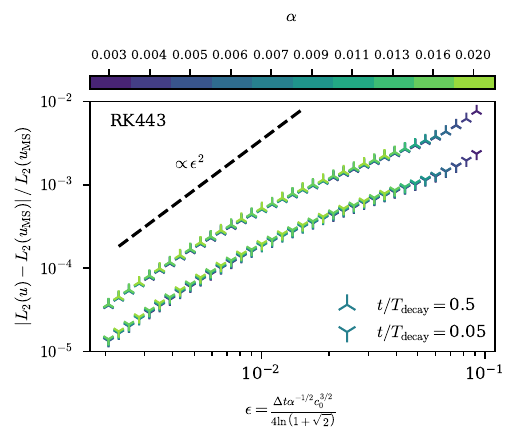}
    \caption{\textbf{Relative errors between measured norms $L_2(u)$ and predicted norms $L_2(u_{\rm MS})$ at intermediate times in RK443 simulations}. 
    Predictions of the $L_2$ norm for RK443 likewise converge in the $\epsilon \rightarrow 0$ limit. Errors are considered at two intermediate times of $t = 0.05 T_{\rm decay}$ and $t = 0.5 T_{\rm decay}$, where $T_{\rm decay}$ is either the simulation decay time or the predicted decay time for $L_2(u)$ and $L_2(u_{\rm MS})$, respectively.
    We again find that errors scale as $\epsilon^2$ for small $\epsilon$, and that the error level is higher at the later time.%. 
    }
    \label{fig:err-wrt-rk443}
\end{figure}

A further capability of the analysis is making predictions about how the full spatial distribution of the numerical solution varies over time. Figures~\ref{fig:err-wrt-group} and~\ref{fig:err-wrt-rk443} demonstrate this capability by comparing predictions and measurements of the $L_2$ norm of the numerical solution at intermediate times. 

In Figure~\ref{fig:err-wrt-group}, we consider the relative error between the predicted norm $L_2(u_{\rm MS})$ at a fraction of the predicted blow-up time $T_{\rm blowup, MS}$ and the measured norm $L_2(u)$ at the same fraction of the actual blow-up time $T_{\rm blowup}$ found in simulations with SBDF1, SBDF2, and RK222. There are two contributions to this error: the difference between $u_{\rm MS}$ and $u$ at the same time---which converges uniformly in time provided the solution has compact support---and the discrepancy between $T_{\rm blowup, MS}$ and $T_{\rm blowup}$. As both of these predictions are expected to have next corrections two orders higher in $\epsilon$ (e.g., as seen for the blow-up times in Figure~\ref{fig:t-ends-rel-err-group}), these two contributions should and do result in similar $\epsilon^2$ convergence rates in the relative error shown here. We still see some discrepancies at small $\epsilon$ in SBDF1, again likely due to the development of asymmetries.

Figure~\ref{fig:err-wrt-rk443} likewise shows the expected $\epsilon^2$ convergence of the $L_2$ relative error for the RK443 timestepping scheme. This result differs from the other three schemes in that here the fractional times are chosen relative to the decay time. Despite this slightly different interpretation, we note that the relative errors in Figure~\ref{fig:err-wrt-rk443} are as low as the $10^{-5}$ to $10^{-4}$ level when $\epsilon \approx 2\times 10^{-3}$, which is similar to how SBDF2 and RK222 perform in that regime. It is unsurprising that the relative error levels are comparable, since the blow-up time scales of SBDF2 and RK222 are within an order of magnitude of the decay time scales for a given choice of $\Delta t$ and $\alpha$. These results demonstrate that our framework is able to predict the progression of different spurious features of the numerical solution over time, which may provide avenues to track and correct for these developments throughout the course of a simulation. 

\section{Conclusions}\label{sec:conc}

This work investigates a pernicious nonlinear timestepping instability found in simulations of solitons. The time scale of the instability is unaffected by the spatial discretization, but depends on the timestepping scheme. This rules out a CFL interpretation of the instability, and we also find that it is not explained by modal predictions from linear stability theory. Instead, we show that a novel multi-scale asymptotic framework is able to recover the evolution of the instability in the semi-discrete limit and predict the blow-up time scale in affected schemes. The asymptotic framework is also able to distinguish between those timestepping schemes which terminate in finite-time blow-up and those which exhibit bounded decay. We attribute the success of our asymptotic analysis to its ability to capture the slow, nonlinear accumulation of timestepping errors, and we find that the von Neumann analysis from linear theory fails because it neglects nonlinearities, which are an important source of error in our simulations even in the $\Delta t \rightarrow 0$ limit.

The asymptotic analysis in Section~\ref{sec:ms} captures all stages of the instability by retaining nonlinear terms and allowing for the numerical solution to evolve on multiple time scales. When evolving soliton solutions of the KdV equation with explicit advection and implicit dispersion, the analysis correctly predicts that the accumulation of timestepping errors will result in changes in the amplitude, speed, and width of the soliton due to predicted time variation of the effective soliton parameter $c$. Our work demonstrates that the order of accuracy of the IMEX scheme sets the slow time scale on which changes in $c$ develop. For some schemes (SBDF1, SBDF2, and RK222), we predict $c$ will approach a finite-time singularity corresponding to soliton blow-up, while in RK443 we predict the soliton will decay. We compare these predictions to numerical results in Section~\ref{sec:results} and confirm that they become more accurate in the $\Delta t \rightarrow 0$ limit. 
 
While this study focuses on the evolution of KdV solitons under a particular set of IMEX timestepping schemes, the multiple scales framework developed here can be readily applied in other contexts. The analysis is already well-equipped for the study of many varieties of multi-step and multi-stage schemes, and it has the potential to be generalized to other classes of timesteppers. For example, Appendix~\ref{sec:ms-im} shows that the analysis is readily applied to fully-implicit multi-step schemes, and predicts that they are more stable than their multi-step IMEX counterparts. In addition to our results on KdV solitons, we report that some of the IMEX schemes in this work also exhibit finite-time blow-up in pseudospectral simulations of a nonlinear generalized Klein-Gordon equation and its cnoidal wave solutions. Much like the solitons in this work, these traveling waves exhibit similar instabilities where the growth in the numerical solution terminates in a power-law singularity at a time sensitive to the timestep size. Future work could confirm that the multiple scales analysis would recover these observations. We have so far found this instability in simulations on periodic intervals where propagating waves have sufficient time to accumulate temporal errors, but similar instabilities may arise in other spatial domains or under different boundary conditions depending on the physical time scales of the problem. The multiple scales framework has been successful when working with initial conditions which exactly solve the continuous problem, but we expect that the method can be extended to the treatment of simulations with arbitrary initial data.

This work demonstrates the shortcomings of linear analysis of timestepping stability in nonlinear simulations where the small scales are set by dispersion and there is no dissipation. If a small amount of dissipation were introduced in the KdV problem, we expect that linear theory would correctly find many of the IMEX schemes considered in this work to be stable in the $\Delta t \rightarrow 0$ limit, as dissipation would eventually overpower the spurious nonlinear growth and therefore cause soliton solutions to decay. Hence, sufficient reduction of the timestep size in problems with a fixed amount of dissipation should be able to avoid nonlinear timestepping instabilities like those considered in this work. It could also be possible to develop a prescription for how much dissipation is necessary to counteract the leading order slow-time instability for a given timestep size. Future work could employ the multiple scales framework to address this idea.

Unexpected behavior in numerical simulations should be met with scrutiny. We bring attention to this intriguing instability to demonstrate the potential shortcomings of linear theory and to incentivize the more careful consideration of the stability of timestepping schemes. We also hope that this work can be referenced as a case study for the utility of asymptotic methods in numerical analysis, as the multiple scales analogy is fundamental to the understanding of this instability.

\section{Acknowledgements}\label{sec:acknowledgements}
The authors thank Emma Kaufman, Liam O'Connor, Erin Cox, Bryan Scott, Caitlin Witt, Whitney Powers, Adrian Fraser, Kyle Auguston, and Ilaria Fontana for helpful discussions and valuable feedback. We also thank the developers of the {\tt Dedalus} project.

B.A.H. is supported by the National Science Foundation Graduate Research Fellowship (Grant no. DGE-2234667). 
D. L. is supported in part by NASA HTMS grant 80NSSC20K1280, NASA OSTFL grant 80NSSC22K1738, Sloan Foundation grant FG-2024-21548, and Simons Foundation grant SFI-MPS-T-MPS-00007353. E.H.A acknowledges partial support from NSF grant PHY-2309135 and Simons Foundation grant (216179, LB) to the Kavli Institute for Theoretical Physics (KITP).
K.J.B. is supported in part by NASA HTMS grant 80NSSC20K1280 and NASA OSTFL grant 80NSSC22K1738.

This research was supported in part through the computational resources and staff contributions provided for the Quest high performance computing facility at Northwestern University which is jointly supported by the Office of the Provost, the Office for Research, and Northwestern University Information Technology.

\section{Data Availability}
Code and processed data used in this work is available in an online public GitHub repository which can be found at \href{https://github.com/benjaminhyatt/nonlin-wave-timestepping-instability}{https://github.com/benjaminhyatt/nonlin-wave-timestepping-instability}.

\appendix

\setcounter{table}{0}
\renewcommand{\thetable}{\Alph{section}\arabic{table}}

\section{IMEX timestepping schemes}\label{sec:IMEX}

We make use of two classes of IMEX timestepping schemes in this work. IMEX schemes evolve the problem
\begin{equation}\label{eq:IMEX-splitting}
    \partial_t u + {\cal L} u = f(u)
\end{equation}
forward in time with implicit evaluation of the linear term ${\cal L} u$ and explicit evaluation of $f(u)$ which may be nonlinear. In general, this splitting is adopted in order to avoid both explicit stability constraints on the oftentimes stiff linear term and possibly expensive iterative methods needed to implicitly evaluate nonlinearities on each timestep. Here we introduce the general structures of the IMEX schemes used in this work. We direct interested readers to the Supplemental Material to see the schemes written out for the KdV equation~\eqref{eq:KdV}, where we take $f(u)=-u\partial_x u$ and ${\cal L} = \alpha \partial_x^3$.

The first class of IMEX schemes we work with is semi-implicit backward differencing (SBDF) schemes. We use multi-step SBDF schemes ranging from first-order to fourth-order accuracy in time (i.e., SBDF1, SBDF2, SBDF3, and SBDF4). These schemes correspond to VSSBDF1-4 as defined in \cite{Wang2008} in the case the timestep size is fixed. An $s$-step SBDF scheme has the form 
\begin{equation}\label{eq:sbdf}
    \frac{1}{\Delta t}\sum_{i=0}^s a_i U^{(n+i)} + \sum_{i=0}^{s} \gamma_i {\cal L} U^{(n+i)} = \sum_{i=0}^{s-1} \beta_i f(U^{(n+i)}),
\end{equation}
where the real coefficients $a_i, \gamma_i$, and $\beta_i$ are chosen to satisfy order conditions such that \eqref{eq:sbdf} is $s$-order accurate in time. The values of the coefficients are given in Table~\ref{tbl:sbdf-coeffs}. We note that $\beta_s$ is not included because a nonzero $\beta_s$ would make the integration of $f$ implicit. A further simplification we can make in this work is to consider purely-implicit integrations of ${\cal L} u$ by taking $\gamma_s = 1$ and $\gamma_i = 0$ for $i = 0, \dots, s-1$.

In order to leverage a multi-step SBDF scheme (i.e.,~\eqref{eq:sbdf} with $s > 1$), one needs to have data for the first $s-1$ iterations. In our simulations we use bootstrapping: taking the first step with SBDF1 and raising the order of accuracy with each subsequent step.
We opt to neglect this bootstrapping procedure in our analysis because our focus is the stability of successive steps of a particular scheme. 

The second class of IMEX schemes we use is mixed explicit, diagonally-implicit Runge-Kutta (ERK+DIRK) schemes. Unlike SBDF schemes which are multi-step in order to obtain higher orders of accuracy, Runge-Kutta (RK) schemes improve in accuracy by being multi-stage. The use of internal stages to determine the solution at the next timestep is sometimes preferred in that no bootstrapping is required. 

A $q$-stage fixed timestep size ERK+DIRK scheme for \eqref{eq:IMEX-splitting} advances by solving stages of the form
\begin{equation}\label{eq:rk}
    \left(1 + \Delta t \tilde{a}_{i,i} {\cal L} \right) U^{(n,i)} = U^{(n,0)} + \Delta t \sum_{j=0}^{i-1} \left(\hat{a}_{i, j} f(U^{(n,j)}, t_n + \tilde{c}_j \Delta t)) - \tilde{a}_{i, j} {\cal L} U^{(n,j)}\right), \quad i = 1, \dots, q.
\end{equation}
We make use of two schemes of this form given in \cite{Ascher1997}, including the two-stage ($q = 2$) second-order scheme RK222 ((2.6) in \cite{Ascher1997}) and the four-stage ($q = 4$) third-order scheme RK443 ((2.8) in \cite{Ascher1997}). In our notation, $\hat{a}_{i,j}$ and $\tilde{a}_{i,j}$ denote the coefficients of the ERK Butcher tableau and (padded) DIRK Butcher tableau, and $\tilde{c}_j$ the corresponding abscissae. The ERK (Left) and DIRK (Right) Butcher tableaus used in this work are
\begin{equation}\label{eq:RK222table}
    \begin{array}
    {c|ccccc}
    0 & 0 \\
    \gamma & \gamma & 0\\
    1 & \delta & 1-\delta & 0 \\
    \hline
    & \delta & 1-\delta & 0
    \end{array} \quad \begin{array}
    {c|ccccc}
    0 & 0   \\
    \gamma & 0 & \gamma\\
    1 & 0 & 1-\gamma & \gamma \\
    \hline
    & 0 & 1-\gamma & \gamma
    \end{array} \quad \text{(RK222)}
\end{equation}
where $\gamma = (2-\sqrt{2})/2$ and $\delta = 1 - 1/(2\gamma)$, and
\begin{equation}\label{eq:RK443table}
    \begin{array}
    {c|ccccc}
    0 & 0 \\
    1/2 & 1/2   & 0\\
    2/3 & 11/18 & 1/18 & 0 \\
    1/2 & 5/6   & -5/6 & 1/2 & 0 \\
    1   & 1/4   & 7/4  & 3/4 & -7/4 & 0\\
    \hline
    & 1/4 & 7/4 & 3/4 & -7/4 & 0
    \end{array} \quad \begin{array}
    {c|ccccc}
    0 & 0 &   \\
    1/2 & 0 & 1/2  & \\
    2/3 & 0 & 1/6  & 1/2  &  \\
    1/2 & 0 & -1/2 & 1/2  & 1/2 &  \\
    1   & 0 & 3/2  & -3/2 & 1/2 & 1/2 \\
    \hline
    & 0 & 3/2 & -3/2 & 1/2 & 1/2
    \end{array} \quad \text{(RK443).}
\end{equation}

Note that for both Runge-Kutta schemes, the abscissae $\tilde{c}_j$ (i.e., the entries in the leftmost column of each tableau) are the same across the ERK and DIRK components, such that the explicit and implicit stages are evaluated simultaneously. We also note that in both schemes the solution of the $q$-th stage is itself the desired update to the solution at the next time level, i.e., $U^{(n+1)} = U^{(n,q)}$. 

\begin{table}[h]
\caption{\textbf{SBDF coefficients}. Here we collect the coefficients of each SBDF scheme (of the form \eqref{eq:sbdf}) considered in this work.}
\label{tbl:sbdf-coeffs}
\begin{tabular*}{\tblwidth}{@{} L | L L L @{}}
\toprule
Scheme & $a_i$ & $\beta_i$ & $\gamma_i$ \\ % Table header row
\midrule \midrule
SBDF1 & $a_1 = 1, a_0 = -1$ & $\beta_0 = 1$ & $\gamma_1 = 1, \gamma_0 = 0$ \\
SBDF2 & $a_2 = 3/2, a_1 = -2, a_0 = 1/2$ & $\beta_1 = 2, \beta_0 = -1$ & $\gamma_2 = 1, \gamma_1 = \gamma_0 = 0$ \\
SBDF3 & $a_3 = 11/6, a_2 = -3, a_1 = 3/2, a_0 = -1/3$ & $\beta_2 = 3, \beta_1 = -3, \beta_0 = 1$ & $\gamma_3 = 1, \gamma_2 = \gamma_1 = \gamma_0 = 0$\\
SBDF4 & $a_4 = 25/12, a_3 = -4, a_2 = 3, a_1 = -4/3, a_0 = 1/4$ & $\beta_3 = 4, \beta_2 = -6, \beta_1 = 4, \beta_0 = -1$ & $\gamma_4 = 1, \gamma_3 = \gamma_2 = \gamma_1 = \gamma_0 = 0$\\
\bottomrule
\end{tabular*}
\end{table}

\section{Linear stability regions assuming a constant background velocity}\label{sec:regions}

\begin{figure}[t]
    \centering
    \includegraphics[width = 0.7\linewidth]{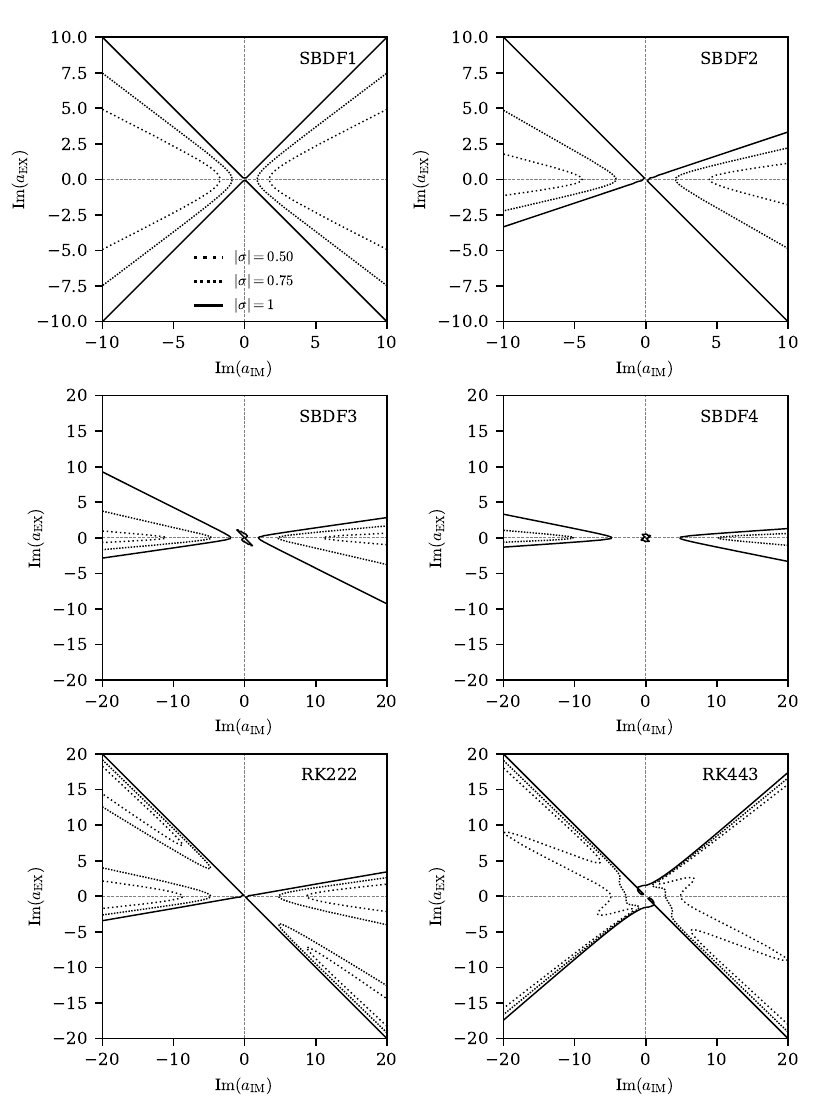}
    \caption{\textbf{Constant-velocity IMEX stability regions}. We plot contours of constant values of the amplification modulus $|\sigma|$ associated with an arbitrary eigenmode of~\eqref{eq:kdv-semi-disc-const}. Positive dispersion ($\alpha > 0$) and advection to the right ($u_0 > 0$) correspond to the bottom-right quadrants of each plot. Large scale (small $k$) modes correspond to values closer to the origin, whereas small scale (large $k$) modes are further away. SBDF3 and SBDF4 have unstable regions near the origin, corresponding to unstable large scale modes. SBDF1, SBDF2, RK222, and RK443 have similar $|\sigma| = 1$  stability boundaries in their bottom-right quadrants. We can infer in each of these schemes that dispersion is stabilizing at small scales (i.e., as $k \rightarrow \infty$), but stability at large scales will depend on the balance between $\alpha$ and $u_0$. RK443 has slightly more stability coverage near the origin than the other schemes.
    }
    \label{fig:IMEX-SR}
\end{figure}

Here we perform a stability analysis of our IMEX timestepping schemes around a constant background solution $u(x,t) = u_0$. This is a more traditional stability calculation than the analysis of perturbations about a soliton presented in Section~\ref{sec:lin-stab}. The advantage of choosing a constant background solution is that the original problem can then be cast as a system of constant coefficient ODEs in time. Consequently, the implicit and explicit parts of the scheme both admit complex exponential eigenfunctions, allowing us to separate out their respective contributions to stability.  

When working with IMEX schemes, the time discretization is applied to ODEs of the standard form
\begin{equation}\label{eq:std}
    \frac{dy}{dt} = \left(a_{\rm IM} + a_{\rm EX}\right)y, \quad a_{\rm IM}, a_{\rm EX} \in \mathbb{C},
\end{equation}
where the implicit and explicit portions are applied to the $a_{\rm IM}y$ and $a_{\rm EX}y$ terms, respectively.

In the case of KdV, we have the linearized constant coefficient equation for the linear perturbation $v$
\begin{equation}\label{eq:kdv-lin-const}
    \partial_t v + \alpha \partial_x^3 v = -u_0\partial_x v,
\end{equation}
which under a complex Fourier basis %(as in~\eqref{eq:expansions}) 
separates into semi-discrete ODEs of the form
\begin{equation}\label{eq:kdv-semi-disc-const}
    \frac{dv_k}{dt}
    = i\left(\alpha\left(\frac{2\pi k}{L}\right)^3 - u_0 \left(\frac{2\pi k}{L}\right)\right)v_k, \quad k = -N/2 + 1, \dots, N/2.%k = 0, \dots, N/2, -N/2 + 1, \dots, -1.
\end{equation}
Assuming implicit evaluation of dispersion and explicit evaluation of advection, we have $a_{\rm IM} = i \alpha\left(\frac{2\pi k}{L}\right)^3$ and $a_{\rm EX} = -i u_0 \left(\frac{2\pi k}{L}\right)$. When we next discretize~\eqref{eq:kdv-semi-disc-const} in time we take the further step of normalizing $a_{\rm IM}$ and $a_{\rm EX}$ by a common factor of $\Delta t$. Since ${\rm Re}(a_{\rm IM}) = {\rm Re}(a_{\rm EX}) = 0$, we can then look at the stability regions in terms of just two dimensions respectively parameterized by the (normalized) imaginary components ${\rm Im}(a_{\rm IM})$ and ${\rm Im}(a_{\rm EX})$.

We plot stability regions for each of the six IMEX schemes considered in this work in Figure~\ref{fig:IMEX-SR}. The stability regions of SBDF3 and SBDF4 are somewhat lacking in coverage near the origin, implying that regardless of timestep size there can be unstable large scale modes. This suggests it is a poor idea to use SBDF3 and SBDF4 to implement implicit dispersion and explicit advection.

The stability regions of the four remaining schemes offer limited insight. They share in common that implicit dispersion should stabilize small scale modes, as for sufficiently large $k$ the positive ${\rm Im}(a_{\rm IM}) \propto k^3$ dominates the negative ${\rm Im}(a_{\rm EX}) \propto k$. Otherwise, there is little that can be ascertained about stability of large and intermediate-scale modes. For smaller values of $k$, modes will be unstable unless $u_0$ is sufficiently small relative to $\alpha$. For instance, if $k = 1$, then $u_0 < \left(2\pi/L\right)^2 \alpha$ is required for stability. That is, when the amount of dispersion set by $\alpha$ is too small, these stability regions can't be accessed by the larger scale modes. If we were to associate $u_0$ with a typical maximum velocity of the soliton (e.g., $u_0 = 3c_0 = 1.5$), this condition would be violated in our simulations based on our selection of $\alpha$ and $L$. However, the results of our multiple scales analysis presented in Section~\ref{sec:results} (e.g., Table~\ref{tbl:ms-results}) suggest that while reduction in $c_0$ may slow down the development of the numerical instability, it would not prevent it altogether. 

While it is simpler to analyze a constant-velocity background state, we consider our work in Section~\ref{sec:lin-stab} to be a more effective use of linear analysis. Both methods predict unstable behavior, but ultimately neither captures the full extent of the behavior in our simulations.

\section{Eigenspectra of the von Neumann analysis
}\label{sec:evp-examples}

\begin{figure}[h]
    \centering
    
    \includegraphics[width = \textwidth]{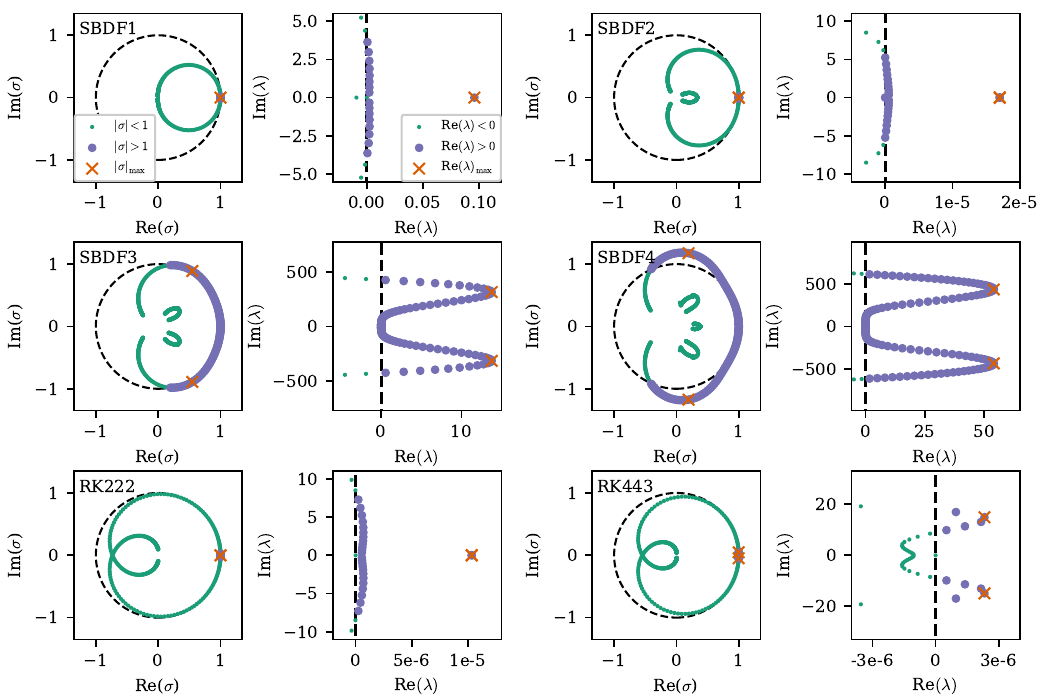}
    \caption{\textbf{von Neumann eigenspectra}. For each scheme, we provide two plots: (Left) spectra of amplification factors $\sigma$ and (Right) a zoomed-in view of the unstable growth rates $\lambda$. We only plot the well-resolved modes from computations with $N=512$. The examples shown correspond to $\alpha = 0.00697$ and $\Delta t = 0.00324$. In SBDF1, SBDF2, RK222, and RK443, unstable amplification factors are clustered near $1$, whereas in SBDF3 and SBDF4 there are overall more unstable modes and they appear at much larger imaginary magnitudes. SBDF1 has a fastest growing mode that is strictly real (i.e., non-oscillatory). SBDF2 and RK222 have complex conjugate pairs of fastest growing modes corresponding to small imaginary parts of about $\pm 10^{-4}$. SBDF3, SBDF4, and RK443 also have complex conjugate pairs with relatively large imaginary parts.}
    \label{fig:evp-dense}
\end{figure}

In Figure~\ref{fig:evp-dense} we provide examples of resolved von Neumann eigenspectra for each IMEX timestepping scheme considered in this work. We determine which eigenmodes are well-resolved by solving the eigenvalue problem with $N_1 = 512$ and $N_2 = 768$ and rejecting amplification factors $\sigma$ that fail to satisfy the ``drift'' condition
\begin{equation}\label{eq:drift}
    \delta_i \equiv \min_{1 \le j \le N_2}\frac{|\sigma_i(N_1) - \sigma_j(N_2)|} {d_i} > 10^{3},
\end{equation}
where the $d_i$ are the intermodal separations in the $N_1$ spectrum, following the approach in Chapter 7 of \cite{Boyd2000}. We find that drift ratios $\delta_i \lesssim 10^{3}$ correspond to poor agreement between the $N_1$ and $N_2$ eigenfunctions, suggesting the former are under-resolved. Between 10\% to 50\% of the modes are found to be under-resolved according to this condition, with higher-order schemes having more rejections. Based on our computations at larger $N$, the missing modes are in general stable and clustered at small $|\sigma|$, so these modes do not impact our arguments about the fastest growing modes in Sections~\ref{sec:lin-stab-results} and~\ref{sec:L-vs-NL}.

Figure~\ref{fig:evp-dense} shows that each scheme has several unstable eigenmodes, corresponding to $|\sigma| > 1$ and ${\rm Re}(\lambda) > 0$. The structure of the eigenspectra of amplification factors $\sigma$ are markedly different in SBDF3 and SBDF4, which exhibit strong linear instabilities in our simulations, compared to the remaining schemes. In SBDF1, SBDF2, RK222, and RK443, the unstable $\sigma$ are each clustered near $1$, %$1 + 0i$,
whereas the spectra of SBDF3 and SBDF4 bulge outside of the unit circle at comparatively large imaginary values.

An additional consideration in our construction of the von Neumann eigenvalue problem in Section~\ref{sec:lin-stab} is how we represent the non-constant background soliton $u_{\rm ex}$ via matrix multiplication with respect to a Fourier basis. To ensure that this representation is valid, we also inspect how truncation of the matrix multiplication affects the determination of the fastest growing eigenmode.

\begin{figure}[t]
    \centering
    \includegraphics[width=\linewidth]{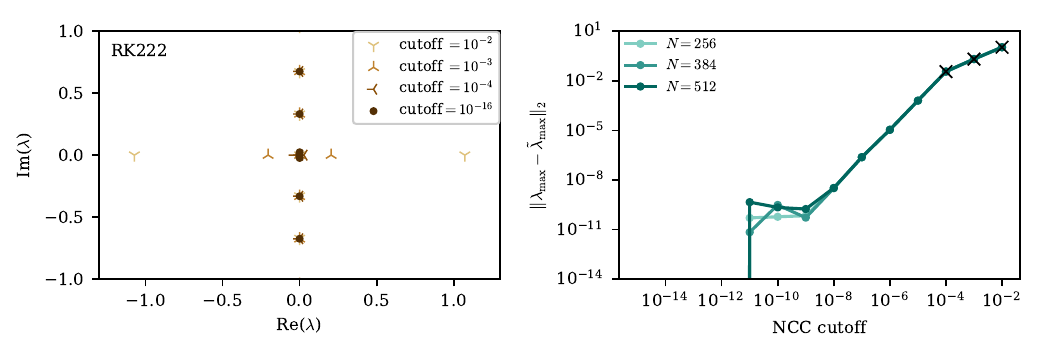}
    \caption{\textbf{Safeguarding against spurious eigenmodes}. 
    (Left) RK222 eigenmodes $\lambda$ for different cutoff values when $N = 512$. (Right) Absolute errors between (a) the $\lambda$ with largest real part when the cutoff is $10^{-16}$, called $\tilde{\lambda}_{\rm max}$ and (b) this value of $\lambda$ for a range of cutoff values, denoted $\lambda_{\rm max}$. When the cutoff value is $10^{-2}, 10^{-3},$ or $10^{-4}$, the identified $\lambda_{\rm max}$ is spurious. When the cutoff value is $10^{-5}$ or smaller, $\lambda_{\rm max}$ is correct and we see fast convergence to $\tilde{\lambda}_{\rm max}$.}
    \label{fig:evp-spurious}
\end{figure}

Taking RK222 for example, Figure~\ref{fig:evp-spurious} plots the eigenmodes $\lambda$ for different ``cutoff'' thresholds which set the amplitude in the space of spectral coefficients to truncate the non-constant background soliton before representing its multiplication as a matrix. When the cutoff value is $\gtrsim 10^{-4}$, there is a spurious pair of real-valued eigenvalues on either side of $\lambda = 0$ that overshadows (i.e., has larger modulus than) the true fastest-growing von Neumann mode. Reduction of the cutoff value to $\lesssim 10^{-5}$ is sufficient to suppress the spurious modes to be smaller than the true fastest growing mode, although cutoff values $\lesssim 10^{-11}$ are needed to resolve the spurious mode to floating-point precision. Figure~\ref{fig:evp-spurious} further shows that the selection of the fastest growing mode is unaffected by changes in $N$, indicating these modes are spatially well-resolved. We conclude that our representation of the non-constant background soliton is valid provided the matrix multiplication is carried out to sufficiently small precision, so our eigenvalue problems in Section~\ref{sec:lin-stab} therefore implement a cutoff value of $10^{-16}$.

\section{Estimates of linear and nonlinear contributions to early error growth}\label{sec:NL-arg}

Section~\ref{sec:lin-stab-results} argues against the validity of the linear stability analysis by comparing estimates for the sizes of linear and nonlinear contributions to growth of numerical error at early times. Here we estimate these two types of contributions. For linear contributions, we consider the error growth rate associated with the fastest growing eigenvalue determined in the linear analysis. For nonlinear contributions, we consider the norm of a nonlinear error term at a fixed early time. The former represents the rate at which perturbations to the background soliton $u_{\rm ex}$ are predicted to grow when nonlinearities are assumed to be negligible. The latter describes the size of the term that is dropped to form the linearized problem~\eqref{eq:KdV-lin}. 

\begin{figure}[h]
    \centering
    \includegraphics{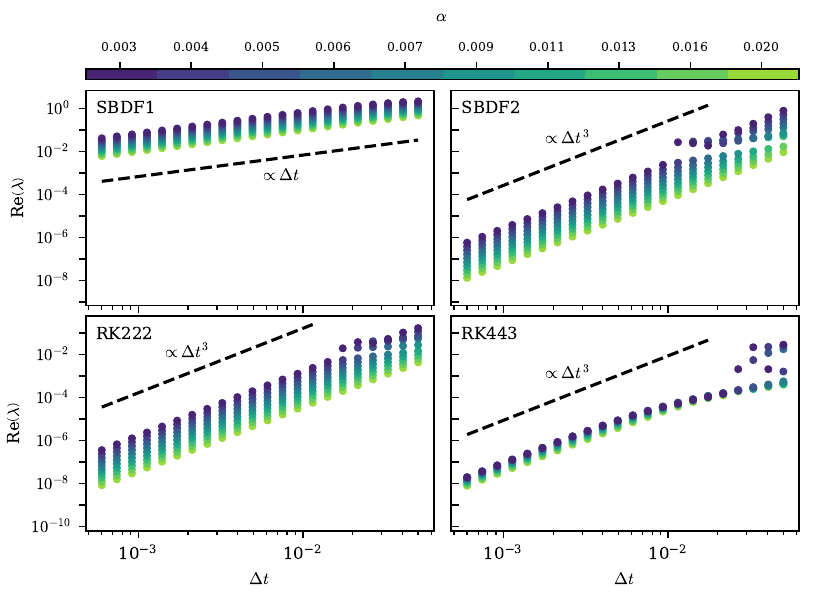}
    \caption{\textbf{Linear instability error growth rates predicted by von Neumann (VN) analysis for SBDF1, SBDF2, RK222, and RK443}. Similar to Figure~\ref{fig:VN-sbdf3-4}, shown are exponential growth rates $\lambda$ inferred from the largest eigenvalues in the von Neumann spectra of SBDF1, SBDF2, RK222, and RK443. These schemes exhibit non-modal growth at early times, growing at paces far exceeding the rates predicted here by the VN analysis. We find ${\rm Re}(\lambda) \sim {\cal O}(\Delta t)$ for SBDF1 and ${\rm Re}(\lambda) \sim {\cal O}(\Delta t^3)$ for SBDF2, RK222, and RK443. 
    }
    \label{fig:VN-sbdf1-2-rk2-3}
\end{figure}

For each of the four timestepping schemes of interest, Figure~\ref{fig:VN-sbdf1-2-rk2-3} gives the real part of the growth rate $\lambda$ associated with the eigenvalue $\sigma$ of largest modulus for a given $(\alpha, \Delta t)$. For each scheme, the dependence of ${\rm Re}(\lambda)$ on $\Delta t$ is given by a power law. A corresponding asymptotic estimate for the change in the norm of the numerical error is formed by considering how ${\rm Re}(\lambda)\|u_{\rm err}\|$ scales with $\Delta t$ at a fixed early time $t_{\rm early}$. 

\begin{figure}[t]
    \centering
    \includegraphics{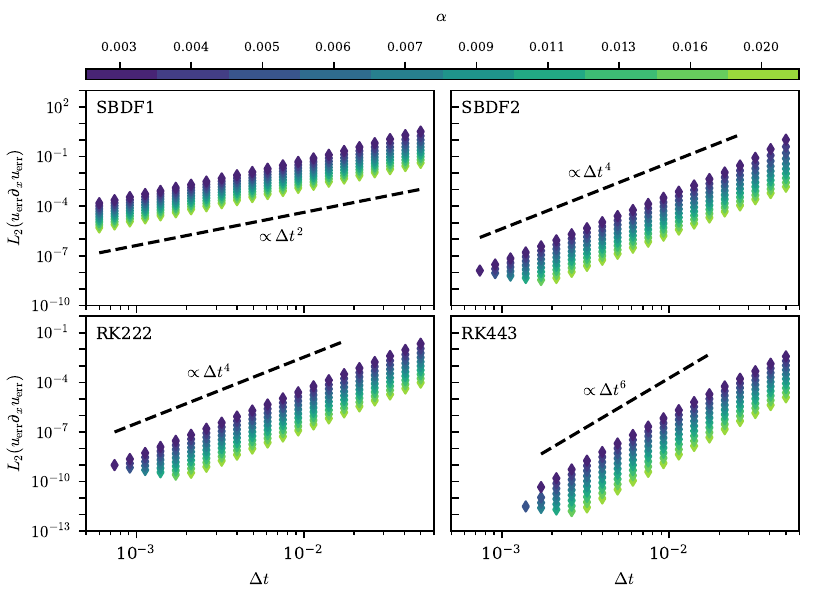}
    \caption{\textbf{Sizes of nonlinear error terms in simulations of solitons}. Shown is the $L_2$ norm of the nonlinear error term $u_{\rm err} \partial_x u_{\rm err}$ measured at a common early time of $t = 0.5$. We find that the error is predominantly determined by the temporal discretization, such that $\|u_{\rm err}\|$ for a fixed $t > 0$ is ${\cal O}(\Delta t^s)$ where $s$ is the order of accuracy of the timestepping scheme. These measurements verify that $\|u_{\rm err} \partial_x u_{\rm err}\|$ is %generally
    ${\cal O}(\Delta t^{2s})$, which is expected as it is quadratic in $u_{\rm err}$.}
    \label{fig:13-NL}
\end{figure}

Figure~\ref{fig:13-NL} reports the norms of $\|u_{\rm err}\partial_x u_{\rm err}\|$ at the same fixed time $t_{\rm early}$. We find that $\|u_{\rm err}\partial_x u_{\rm err}\|$ likewise has a power law scaling in $\Delta t$, and is in fact consistent with the square of the scaling in $\Delta t$ predicted for $\|u_{\rm err}\|$ based on the order of accuracy of the timestepping scheme. Therefore, the comparison between the two contributions ${\rm Re}(\lambda)\|u_{\rm err}\|$ and $\|u_{\rm err}\partial_x u_{\rm err}\|$ ultimately reduces to that between ${\rm Re}(\lambda)$ and $\|u_{\rm err}\|$ at a common early time. Put together, we summarize the comparison of these two contributions in Table~\ref{tbl:err-L-vs-NL}. Despite being solely based on the scaling of the contributions with $\Delta t$ (i.e., without considering %knowing
prefactors), we find that for the problem at hand this argument provides a consistent explanation for which schemes develop a timestepping instability incompatible with the %modal
predictions of linear theory. 

\section{Predictions of the multiple scales analysis for fully-implicit schemes}\label{sec:ms-im}
This study focuses on the stability of IMEX schemes because they are a widely used class of methods and they are where we first observed the timestepping instability of interest. Here we briefly consider the extension of the multiple scales analysis to a few fully-implicit schemes. For simplicity, we consider the fully-implicit analogs of SBDF1 and SBDF2, which are BDF1 (Backward Euler) and BDF2. Written for the KdV equation~\eqref{eq:KdV}, these are given by
\begin{equation}\label{eq:bdf1-kdv}
    \frac{1}{\Delta t}\left(U^{(n+1)}-U^{(n)}\right) + \alpha \partial_x^3 U^{(n+1)} + U^{(n+1)} \partial_x U^{(n+1)} = 0, \quad \text{(BDF1)}
\end{equation}
and
\begin{equation}\label{eq:bdf2-kdv}
    \frac{1}{\Delta t}\left(\frac{3}{2}U^{(n+2)}-2U^{(n+1)}+\frac{1}{2}U^{(n)}\right) + \alpha \partial_x^3 U^{(n+2)} + U^{(n+2)} \partial_x U^{(n+2)} = 0, \quad \text{(BDF2)}
\end{equation}
respectively. As~\eqref{eq:bdf1-kdv} and~\eqref{eq:bdf2-kdv} now treat the nonlinear advection term implicitly, we would need to adopt a suitable iterative method in order to do numerical computations with either of these schemes. This study does not include computations with these schemes, but we can use the multiple scales approach to predict whether or not these fully-implicit schemes are unstable to the nonlinear timestepping instability. 

Following the same procedure as outlined in Sections~\ref{sec:ms} and~\ref{sec:results}, we find that BDF1 has the expected fast-time behavior given in~\eqref{eq:ms-leading-order}, and at ${\cal O}(\epsilon)$, BDF1 admits a nontrivial solvability condition of the form 
\begin{equation}\label{eq:bdf1-sc}
    \frac{1}{2}\partial_{t_1}{\left(\int_{-L/2}^{L/2} U_0^2\,d\xi\right)} = \int_{-L/2}^{L/2} U_0\tau_0\partial_{t_0}{\left(-\frac{\alpha}{2}\partial_x^3U_0 - \frac{1}{2}U_0\partial_xU_0\right)}\, d\xi.
\end{equation}
When we assume $U_0$ has the form of the soliton given in~\eqref{eq:ms-U0-general} and consider the $L \rightarrow \infty$ limit, we find that~\eqref{eq:bdf1-sc} corresponds to the separable ODE
\begin{equation}\label{eq:bdf1-ode-related}
    \frac{\partial c}{\partial t_1} = -\frac{2}{15}\frac{\tau_0}{\alpha}c^4,
\end{equation}
which with $c(0)=c_0$ admits the power law solution
\begin{equation}\label{eq:bdf1-PL-related}
    c(t_1) = \left(c_0^{-3} + \frac{6}{15}\frac{\tau_0}{\alpha}t_1\right)^{-1/3} \quad \Leftrightarrow \quad c(t) = \left(c_0^{-3} + \frac{6}{15}\frac{\Delta t}{\alpha}t\right)^{-1/3}.
\end{equation}
Similarly, BDF2 has the same form as SBDF2 at orders ${\cal O}(1)$ and ${\cal O}(\epsilon)$, and at ${\cal O}(\epsilon^2)$ has a similar trivial solvability condition due to parity. Finally, at ${\cal O}(\epsilon^3)$, BDF2 has the nontrivial solvability condition
\begin{equation}\label{eq:bdf2-e3-sc}
    \frac{1}{2}\partial_{t_3}{\left(\int_{-L/2}^{L/2}U_0^2\,d\xi\right)} = \int_{-L/2}^{L/2} U_0\tau_0^3\partial_{t_0}^3{\left(\frac{\alpha}{4}\partial_x^3U_0 + \frac{1}{4}U_0\partial_xU_0\right)}\, d\xi,
\end{equation}
which in the $L\rightarrow \infty$ limit corresponds to the separable ODE
\begin{equation}\label{eq:bdf2-ode-related}
    \frac{\partial c}{\partial t_3} = -\frac{1}{21}\frac{\tau_0^3}{\alpha^2}c^7,
\end{equation}
which with $c(0)=c_0$ admits the power law solution
\begin{equation}\label{eq:bdf2-PL-related}
    c(t_3) = \left(c_0^{-6} +\frac{2}{7}\frac{\tau_0^3}{\alpha^2}t_3\right)^{-1/6} \quad \Leftrightarrow \quad c(t) = \left(c_0^{-6} +\frac{2}{7}\frac{\Delta t^3}{\alpha^2}t\right)^{-1/6}.
\end{equation}

\eqref{eq:bdf1-PL-related} and~\eqref{eq:bdf2-PL-related} respectively indicate that the multiple scales analysis predicts BDF1 and BDF2 are stable to the nonlinear timestepping instability. While this is not definitive proof of their stability, we can gain some useful insights from this exercise.

We note that the only difference between~\eqref{eq:bdf1-sc} and~\eqref{eq:bdf2-e3-sc} and their IMEX counterparts is that the coefficient of the advection term $U_0\partial_x U_0$ is now the same as that of the dispersion term $\alpha\partial_x^3 U_0$, whereas in~\eqref{eq:sbdf1-sc} and~\eqref{eq:sbdf2-e3-sc} the coefficients differ. This is a natural consequence of both terms entering the fully-implicit scheme in the same way, i.e., they both enter according to $\gamma_i$ in~\eqref{eq:sbdf}.

Furthermore, comparing the power law results of the fully-implicit scheme BDF2 with that of the stable IMEX scheme RK443, we see that the prefactor in front of time $t$ in~\eqref{eq:bdf2-PL-related} of $(2/7)\Delta t^3 \alpha^{-2}$ is approximately $9$ times smaller than the prefactor $(77069/30030)\Delta t^3 \alpha^{-2}$ in~\eqref{eq:rk443-PL-related}, indicating that the slow-time evolution of the soliton under BDF2 would proceed at a slower pace by nearly an order of magnitude.

These observations lead us to speculate that when the advection term is treated implicitly alongside the dispersion term, the slow-time evolution of the soliton is stable, but when the advection term is treated explicitly, there is competition between advection and dispersion which accelerates the slow-time evolution and sometimes leads to instability. At this time, this interpretation is a conjecture, but this exercise suggests that it would be worthwhile to follow-up on this work to better understand why schemes such as RK443 (and likely BDF1 and BDF2) are stable.

\printcredits

%% Loading bibliography style file
%\bibliographystyle{model1-num-names}
\bibliographystyle{cas-model2-names}

% Loading bibliography database
\bibliography{cas-refs}

\end{document}

\title{Supplemental Material for \\
``Multiple scales analysis of a nonlinear timestepping instability\\in simulations of solitons"}
\author{Benjamin A. Hyatt, Daniel Lecoanet, Evan H. Anders, Keaton J. Burns}
\date{}
\maketitle

\tableofcontents

\section{Implicit-explicit (IMEX) timestepping schemes used in this work}

We use IMEX timestepping schemes to evolve problems of the form
\begin{equation}\label{eq:IMEX-splitting}
    \partial_t u + {\cal L} u = f(u)
\end{equation}
by treating ${\cal L} u$ implicitly and $f(u)$ explicitly. In particular, when \eqref{eq:IMEX-splitting} is the nonlinear KdV equation
\begin{equation}\label{eq:kdv-sm}
    \partial_t u + \alpha \partial_x^3 u = -u\partial_x u, 
\end{equation}
we have ${\cal L} u = \alpha \partial_x^3 u$ and $f(u) = -u \partial_x u$. Here we articulate how two classes of IMEX schemes---semi-implicit backward differencing (SBDF) schemes and mixed explicit, diagonally-implicit Runge-Kutta (ERK+DIRK) schemes---are used to discretize the KdV equation \eqref{eq:kdv-sm}. 

We denote different SBDF schemes as SBDF$(s)$ where $s$ is the order of accuracy. This follows~\cite{Wang2008} in the case of a fixed timestep size $\Delta t$. Following~\cite{Ascher1997}, we denote the ERK+DIRK schemes as RK$(q,r,s)$ where $q$ and $r$ are the number of implicit and explicit stages to be solved respectively, and $s$ is again the order of accuracy. 

\subsection{Semi-implicit backward difference schemes (SBDF)}

SBDF$(s)$ corresponds to an $s$-step SBDF scheme for \eqref{eq:IMEX-splitting} which has the form 
\begin{equation}\label{eq:sbdf}
    \frac{1}{\Delta t}\sum_{i=0}^s a_i U^{(n+i)} + \sum_{i=0}^{s} \gamma_i {\cal L} U^{(n+i)} = \sum_{i=0}^{s-1} \beta_i f(U^{(n+i)}),
\end{equation}
with the fixed timestep size $\Delta t$. This is $s$-order accurate provided the coefficients $a_i, \beta_i,$ and $\gamma_i$ satisfy the necessary order conditions (see (2.4) in \cite{Wang2008}). 

SBDF1, otherwise known as Forward-Backward Euler, corresponds to scheme (2.6) in \cite{Wang2008}. Its coefficients are
\begin{equation}\label{eq:sbdf1-coeffs}
    a_1 = 1, a_0 = -1, \quad \beta_0 = 1, \quad \gamma_1 = 1, \gamma_0 = 0,
\end{equation}
which leads to the time discretization of the KdV equation
\begin{equation}\label{eq:sbdf1-kdv}
    \frac{1}{\Delta t}\left(U^{(n+1)}-U^{(n)}\right) + \alpha \partial_x^3 U^{(n+1)} = -U^{(n)} \partial_x U^{(n)}.
\end{equation}

SBDF2 corresponds to scheme (2.8) in \cite{Wang2008} in the case of constant step size. Its coefficients are
\begin{equation}\label{eq:sbdf2-coeffs}
    a_2 = 3/2, a_1 = -2, a_0 = 1/2, \quad \beta_1 = 2, \beta_0 = -1 \quad \gamma_2 = 1, \gamma_1 = \gamma_0 = 0,
\end{equation}
leading to the KdV discretization
\begin{equation}\label{eq:sbdf2-kdv}
    \frac{1}{\Delta t}\left(\frac{3}{2}U^{(n+2)}-2U^{(n+1)}+\frac{1}{2}U^{(n)}\right) + \alpha \partial_x^3 U^{(n+2)} = -2U^{(n+1)} \partial_x U^{(n+1)} + U^{(n)} \partial_x U^{(n)}.
\end{equation}

SBDF3 corresponds to scheme (2.14) in \cite{Wang2008} in the case of constant step size. Its coefficients are
\begin{equation}\label{eq:sbdf3-coeffs}
    a_3 = 11/6, a_2 = -3, a_1 = 3/2, a_0 = -1/3, \quad \beta_2 = 3, \beta_1 = -3, \beta_0 = 1, \quad \gamma_3 = 1, \gamma_2 = \gamma_1 = \gamma_0 = 0,
\end{equation}
leading to the KdV discretization
\begin{equation}\label{eq:sbdf3-kdv}
    \frac{1}{\Delta t}\left(\frac{11}{6} U^{(n+3)} - 3 U^{(n+2)} + \frac{3}{2} U^{(n+1)} - \frac{1}{3} U^{(n)}\right) + \alpha \partial_x^3 U^{(n+3)} = -3 U^{(n+2)}\partial_x U^{(n+2)} + 3 U^{(n+1)}\partial_x U^{(n+1)} - U^{(n)}\partial_x U^{(n)}.
\end{equation}

SBDF4, the last multi-step scheme we consider, corresponds to scheme (2.15) in \cite{Wang2008} in the case of constant step size. Its coefficients are
\begin{equation}\label{eq:sbdf4-coeffs}
    a_4 = 25/12, a_3 = -4, a_2 = 3, a_1 = -4/3, a_0 = 1/4, \quad \beta_3 = 4, \beta_2 = -6, \beta_1 = 4, \beta_0 = -1, \quad \gamma_4 = 1, \gamma_3 = \gamma_2 = \gamma_1 = \gamma_0 = 0,
\end{equation}
leading to the KdV discretization
\begin{equation}\label{eq:sbdf4-kdv}
\begin{aligned}
    \frac{1}{\Delta t}\left(\frac{25}{12}U^{(n+4)} - 4 U^{(n+3)} + 3 U^{(n+2)} - \frac{4}{3} U^{(n+1)} + \frac{1}{4} U^{(n)}\right) + \alpha \partial_x^3 U^{(n+4)} = &-4 U^{(n+3)}\partial_x U^{(n+3)} + 6 U^{(n+2)}\partial_x U^{(n+2)} \\
    &- 4U^{(n+1)}\partial_x U^{(n+1)} + U^{(n)}\partial_x U^{(n)}.
\end{aligned}
\end{equation}

\subsection{Mixed explicit, diagonally-implicit Runge-Kutta schemes (ERK + DIRK)}

RK$(q,r,s)$ corresponds to a $q$-stage ERK+DIRK scheme for \eqref{eq:IMEX-splitting} which advances by solving stages of the form
\begin{equation}\label{eq:rk}
    \left(1 + \Delta t \tilde{a}_{i,i} {\cal L} \right) U^{(n,i)} = U^{(n,0)} + \Delta t \sum_{j=0}^{i-1} \left(\hat{a}_{i, j} f(U^{(n,j)}, t_n + \tilde{c}_j \Delta t) - \tilde{a}_{i, j} {\cal L} U^{(n,j)}\right), \quad i = 1, \dots, q.
\end{equation}
with the fixed timestep size $\Delta t$. The coefficients $\hat{a}_{i,j}$, $\tilde{a}_{i,j}$ are elements of the ERK and (padded) DIRK Butcher tableau, and $\tilde{c}_j$ the corresponding abscissae. In the schemes we work with, the abscissae are the same for the explicit and implicit parts, and hence we solve the $q = r$ stages simultaneously. The scheme is $s$-order accurate provided the coefficients follow the order conditions prescribed in Part II, Chapter 15 of~\cite{Hairer1993}.

RK222 corresponds to scheme (2.6) in \cite{Ascher1997} which gives the corresponding Butcher tableaus, which we also provide in Appendix A of the main paper. This leads to the KdV discretization
\begin{equation}\label{eq:rk222-kdv}
    \begin{aligned}
    \left(1+ \Delta t \gamma \alpha \partial_x^3\right)U^{(n,1)} &= U^{(n,0)} + \Delta t \gamma \left(-U^{(n,0)}\partial_x U^{(n,0)}\right) ,\\
    \left(1+\Delta t \gamma \alpha \partial_x^3\right)U^{(n,2)} &= U^{(n,0)} + \Delta t \delta \left(-U^{(n,0)}\partial_x U^{(n,0)}\right) - \Delta t (1 - \gamma) \alpha \partial_x^3 U^{(n,1)} + \Delta t (1-\delta) \left(-U^{(n,1)}\partial_x U^{(n,1)}\right),
    \end{aligned}
\end{equation}
where $\gamma = (2-\sqrt{2})/2$ and $\delta = 1 - 1/(2\gamma)$

RK443 corresponds to scheme (2.8) in \cite{Ascher1997} which gives the corresponding Butcher tableaus, which can again be found in Appendix A of the main paper. This leads to the KdV discretization
\begin{equation}\label{eq:rk443-kdv}
    \begin{aligned}
    \left(1+ \Delta t \left(\frac{1}{2}\right)\alpha \partial_x^3\right)U^{(n,1)} &= U^{(n,0)} + \Delta t \left(\frac{1}{2}\right)\left(-U^{(n,0)}\partial_x U^{(n,0)}\right) ,\\
    \left(1 + \Delta t \left(\frac{1}{2}\right) \alpha \partial_x^3\right)U^{(n,2)} &= U^{(n,0)} + \Delta t \left(\left(\frac{11}{18}\right)\left(-U^{(n,0)}\partial_x U^{(n,0)}\right) + \left(\frac{1}{18}\right)\left(-U^{(n,1)}\partial_x U^{(n,1)}\right)\right) - \Delta t \left(\frac{1}{6}\right)\alpha \partial_x^3 U^{(n,1)},\\
    \left(1 + \Delta t \left(\frac{1}{2}\right)\alpha \partial_x^3\right)U^{(n,3)} &= U^{(n,0)} + \Delta t \left(\left(\frac{5}{6}\right)\left(-U^{(n,0)}\partial_x U^{(n,0)}\right) + \left(\frac{-5}{6}\right)\left(-U^{(n,1)}\partial_x U^{(n,1)}\right) + \left(\frac{1}{2}\right)\left(-U^{(n,2)}\partial_x U^{(n,2)}\right)\right)\\ 
    &- \Delta t\left(\left(\frac{-1}{2}\right)\alpha \partial_x^3 U^{(n,1)} + \left(\frac{1}{2}\right)\alpha\partial_x^3 U^{(n,2)}\right),\\
    \left(1 + \Delta t \left(\frac{1}{2}\right)\alpha \partial_x^3\right)U^{(n,4)} &= U^{(n,0)} + \Delta t \left(\left(\frac{1}{4}\right)\left(-U^{(n,0)}\partial_x U^{(n,0)}\right) + \left(\frac{7}{4}\right)\left(-U^{(n,1)}\partial_x U^{(n,1)}\right) + \left(\frac{3}{4}\right)\left(-U^{(n,2)}\partial_x U^{(n,2)}\right)\right.\\
    &\left. + \left(\frac{-7}{4}\right)\left(-U^{(n,3)}\partial_x U^{(n,3)}\right)\right) - \Delta t\left(\left(\frac{3}{2}\right)\alpha\partial_x^3 U^{(n,1)} + \left(\frac{-3}{2}\right)\alpha\partial_x^3 U^{(n,2)} + \left(\frac{1}{2}\right)\alpha\partial_x^3 U^{(n,3)}\right).
    \end{aligned}
\end{equation}

\section{Multiple scales analysis equations at each order in $\epsilon = \Delta t / \tau_0$}\label{sec:ms-each}

We study how SBDF1, SBDF2, RK222, and RK443 evolve soliton solutions of \eqref{eq:kdv-sm} by assuming the numerical solution possibly varies on the multiple time scales
\begin{equation}\label{eq:times-sm}
    t_0 = t, \quad t_1 = \epsilon t, \quad \dots, \quad t_s = \epsilon^s t, \quad t_{s+1} = \epsilon^{s+1} t,
\end{equation}
where the small parameter $\epsilon = \Delta t / \tau_0$ is a ratio between the timestep size and soliton propagation time scale. We assume the numerical solution at time $t_n = n \Delta t$ has the regular expansion
\begin{equation}\label{eq:reg-sm}
    U^{(n)}(x) = \sum_{k=0}^{s+1}\epsilon^k U_k^{(n)}(x, t_{0}, \dots, t_{s+1}) + {\cal O}(\epsilon^{s+2}).
\end{equation}
For each scheme we solve for $U_0^{(n)}$, a function of $x$ and some subset of the introduced time scales, as a leading order approximation to $U^{(n)}(x)$. 

In order to solve for $U^{(n)}(x)$, we represent the numerical solution at each discrete time level about time level $t_n$ by way of a Taylor expansion in time. This is a relevant consideration for both multi-step and multi-stage schemes, which will relate the numerical solution for at least more than one time level (e.g., $U^{(n)}(x)$ and $U^{(n+1)}(x)$). For an $s$-step multi-step schemes, we form a Taylor expansion for each of $U^{(n+1)}$, $\dots$, and $U^{(n+s)}$ about $U^{(n)}$. In multi-stage schemes, we need only perform a Taylor expansion on $U^{(n+1)}(x)$ (i.e., the final stage $U_q^{(n)}(x)$). For the internal stages $U_i^{(n)}$ for $i = 1, \dots, q-1$, rather than perform Taylor expansions about $U^{(n)}(x)$ we assume the $i$-th internal stage has a regular expansion of the form
\begin{equation}\label{eq:reg-internal-sm}
    U_i^{(n)}(x) = \sum_{k=0}^{s+1}\epsilon^k U_{i,k}^{(n)}(x, t_{0}, \dots, t_{s+1}) + {\cal O}(\epsilon^{s+2}).
\end{equation}

In the subsections that follow, we report what this procedure yields for each timestepping scheme at each order in $\epsilon$. After reporting these equations, we note any relevant assumptions that were used to drop insignificant terms throughout the solution process. Each assumption is informed by a corresponding solvability condition (as described in the ``Multiple scales analysis'' section of the main paper). Finally, we emphasize which equation was used to obtain a nontrivial solvability condition. We explicitly report these conditions in Section~\ref{sec:sc}. When we assume a finite domain, we obtain solvability conditions by integrating an equation over $\xi \in (-L/2, L/2)$ assuming periodic boundary conditions (where $\xi = x - c t_0$). When we take $L \rightarrow \infty$, we instead integrate over $\xi \in (-\infty, \infty)$ and assume decay at infinity. 

\subsection{SBDF1}
For SBDF1, the multiple scales analysis yields
\begin{subequations}
\begin{align}
    {\cal O}(\epsilon^0) \quad &\partial_{t_0} U_0 + \alpha \partial_x^3 U_0 + U_0 \partial_x U_0 = 0, \label{eq:sbdf1-e0} \\
    {\cal O}(\epsilon^1) \quad &\partial_{t_0} U_1 + \alpha\partial_x^3 U_1 + \partial_x(U_0 U_1) = - \partial_{t_1} U_0 + \tau_0\partial_{t_0}\left(-\frac{\alpha}{2}\partial_x^3U_0 + \frac{1}{2}U_0\partial_xU_0\right). \label{eq:sbdf1-e1}
\end{align}
\end{subequations}
No extra assumptions had to be adopted. We obtain the nontrivial solvability condition by integrating \eqref{eq:sbdf1-e1} with respect to $\xi = x - ct_0$. 

\subsection{SBDF2}
For SBDF2, the multiple scales analysis yields
\begin{subequations}
\begin{align}
    {\cal O}(\epsilon^0) \quad &\partial_{t_0} U_0 + \alpha \partial_x^3 U_0 + U_0 \partial_x U_0 = 0, \label{eq:sbdf2-e0} \\
    {\cal O}(\epsilon^1) \quad &\partial_{t_0} U_1 + \alpha \partial_x^3 U_1 + \partial_x \left(U_0 U_1\right) = -\partial_{t_1} U_0, \label{eq:sbdf2-e1}\\
    {\cal O}(\epsilon^2) \quad &\partial_{t_0} U_2 + \alpha \partial_x^3 U_2 + \partial_x \left(U_0 U_2\right) = -\partial_{t_2} U_0 + \tau_0^2 \partial_{t_0}^2\left(-\frac{\alpha}{3}\partial_x^3 U_0 + \frac{2}{3}U_0\partial_x U_0\right), \label{eq:sbdf2-e2} \\
    {\cal O}(\epsilon^3) \quad &\partial_{t_0}U_3 + \alpha \partial_x^3 U_3 + \partial_x (U_0 U_3) = - \partial_{t_3}U_0 + \tau_0^3 \partial_{t_0}^3 \left(\frac{\alpha}{4} \partial_x^3 U_0 - \frac{3}{4} U_0 \partial_x U_0 \right). \label{eq:sbdf2-e3}
\end{align}
\end{subequations}
For each equation after \eqref{eq:sbdf2-e1}, we adopt $\partial_{t_1} U_0 = 0$ and $U_1 = 0$. After $\eqref{eq:sbdf2-e2}$ we further adopt $\partial_{t_2} U_0 = 0$ and $\partial_{t_1} U_2 = 0$. We obtain the nontrivial solvability condition by integrating \eqref{eq:sbdf2-e3} with respect to $\xi = x - ct_0$.

\subsection{RK222}

For RK222, the multiple scales analysis yields
\allowdisplaybreaks
\begin{subequations}
\begin{align}
    {\cal O}(\epsilon^{-1}) \quad &U_0^{(1)} = U_0, \label{eq:rk222-i1-en1} \\
    {\cal O}(\epsilon^0) \quad &\partial_{t_0} U_0 + \alpha \partial_x^3 U_0 + U_0 \partial_x U_0 = 0, \label{eq:rk222-i2-e0} \\
    &U_1^{(1)} = U_1 - \tau_0 \gamma U_0 \partial_x U_0 - \tau_0 \alpha \gamma \partial_x^3 U_0^{(1)}, \label{eq:rk222-i1-e0} \\
    {\cal O}(\epsilon^1) \quad &\partial_{t_0}U_1 + \alpha \partial_x^3 U_1 + \partial_x(U_0 U_1) = -\partial_{t_1}U_0, \label{eq:rk222-i2-e1}\\
    &U_2^{(1)} = U_2 -\tau_0 \gamma U_1 \partial_x U_0 -\tau_0 \gamma U_0 \partial_x U_1 - \tau_0 \alpha \gamma \partial_x^3 U_1^{(1)}, \label{eq:rk222-i1-e1} \\
    {\cal O}(\epsilon^2) \quad &\partial_{t_0}U_2 + \alpha \partial_x^3 U_2 + \partial_x (U_0 U_2) = - \partial_{t_2}U_0 -\frac{1}{6}\tau_0^2\partial_{t_0}^3 U_0 -\frac{1}{2} \tau_0^2 \gamma U_0(\partial_x U_0)^3 -\frac{1}{2} \tau_0^2 \gamma U_0^2 \partial_x U_0 \partial_x^2 U_0 \label{eq:rk222-i2-e2} \\
    &-\frac{3}{2} \tau_0^2 \alpha \gamma \partial_x U_0 (\partial_x^2 U_0)^2 - \frac{5}{2} \tau_0^2 \alpha \gamma (\partial_x U_0)^2 \partial_x^3 U_0 - \frac{11}{2} \tau_0^2 \alpha \gamma U_0 \partial_x^2 U_0 \partial_x^3 U_0 -\frac{1}{2} \tau_0^2 \alpha \gamma \partial_x^3\partial_{t_0}^2 U_0 \nonumber \\
    &-\frac{7}{2}\tau_0^2 \alpha \gamma U_0 \partial_x U_0 \partial_x^4 U_0 -\frac{1}{2} \tau_0^2 \alpha^2 \gamma \partial_x^3 U_0 \partial_x^4 U_0 -35 \tau_0^2 \alpha^2 \gamma^2 \partial_x^3 U_0 \partial_x^4 U_0 + 35 \tau_0^2 \alpha^2 \gamma^3 \partial_x^3 U_0 \partial_x^4 U_0 \nonumber \\
    &-\frac{1}{2} \tau_0^2 \alpha \gamma U_0^2 \partial_x^5 U_0  -21 \tau_0^2 \alpha^2 \gamma^2 \partial_x^2 U_0 \partial_x^5 U_0 + 21 \tau_0^2 \alpha^2 \gamma^3 \partial_x^2 U_0 \partial_x^5 U_0 - \frac{1}{2} \tau_0^2 \alpha^2 \gamma \partial_x U_0 \partial_x^6 U_0 \nonumber \\
    &-7 \tau_0^2 \alpha^2 \gamma^2 \partial_x U_0 \partial_x^6 U_0 + 7 \tau_0^2 \alpha^2 \gamma^3 \partial_x U_0 \partial_x^6 U_0 -\frac{1}{2} \tau_0^2 \alpha^2 \gamma U_0 \partial_x^7 U_0 - \tau_0^2 \alpha^2 \gamma^2 U_0 \partial_x^7 U_0 \nonumber \\
    &+\tau_0^2 \alpha^2 \gamma^3 U_0 \partial_x^7 U_0 - \tau_0^2 \alpha^3 \gamma^2 \partial_x^9 U_0 + \tau_0^2 \alpha^3 \gamma^3 \partial_x^9 U_0, \nonumber \\ 
    &U_3^{(1)} = U_3 - \tau_0 \gamma U_2 \partial_x U_0 - \tau_0 \gamma U_1 \partial_x U_1  
    -\tau_0 \gamma U_0 \partial_x U_2 -\tau_0 \alpha \gamma \partial_x^3 U_2^{(1)}, \label{eq:rk222-i1-e2} \\ %|_{U_1 = 0}, \label{eq:rk222-i1-e2} \\
    {\cal O}(\epsilon^3) \quad &\partial_{t_0} U_3 + \alpha \partial_x^3 U_3 + \partial_x (U_0 U_3) = - \partial_{t_3} U_0 -\frac{1}{2} \tau_0 \partial_{t_0}^2 U_2 -\frac{1}{24} \tau_0^3 \partial_{t_0}^4 U_0 + \tau_0 U_2 (\partial_x U_0)^2 + 2 \tau_0 U_0 \partial_x U_0 \partial_x U_2 \label{eq:rk222-i2-e3}\\
    &+ \tau_0 U_0 U_2 \partial_x^2 U_0 + \frac{3}{2} \tau_0^3 \alpha \gamma^2 (\partial_x U_0)^2 (\partial_x^2 U_0)^2 +\frac{3}{2} \tau_0^3 \alpha \gamma^2 U_0 (\partial_x^2 U_0)^3 + \frac{1}{2} \tau_0 U_0^2 \partial_x^2 U_2 \nonumber \\
    & + 6 \tau_0 \alpha \gamma \partial_x^2 U_0 \partial_x^2 U_2 -6 \tau_0 \alpha \gamma^2 \partial_x^2 U_0 \partial_x^2 U_2 + 2 \tau_0^3 \alpha \gamma^2 (\partial_x U_0)^3 \partial_x^3 U_0 + \frac{1}{2} \tau_0 \alpha \partial_x U_2 \partial_x^3 U_0\nonumber \\
    &  + 4 \tau_0 \alpha \gamma \partial_x U_2 \partial_x^3 U_0 - 4 \tau_0 \alpha  \gamma^2 \partial_x U_2 \partial_x^3 U_0 + 7 \tau_0^3 \alpha \gamma^2 U_0 \partial_x U_0 \partial_x^2 U_0 \partial_x^3 U_0 +5 \tau_0^3 \alpha^2 \gamma^2 \partial_x^2 U_0 (\partial_x^3 U_0)^2 \nonumber \\
    & +\frac{1}{2} \tau_0 \alpha \partial_x U_0 \partial_x^3 U_2 + 4 \tau_0 \alpha  \gamma \partial_x U_0 \partial_x^3 U_2 -4 \tau_0 \alpha \gamma^2 \partial_x U_0 \partial_x^3 U_2 -\tau_0 \alpha \gamma \partial_{t_0}\partial_x^3 U_2 \nonumber \\
    & -\frac{1}{6} \tau_0^3 \alpha \gamma  \partial_{t_0}^3 \partial_x^3 U_0 +\frac{1}{2} \tau_0 \alpha U_2 \partial_x^4 U_0 + \tau_0 \alpha \gamma  U_2 \partial_x^4 U_0 -\tau_0 \alpha \gamma^2 U_2 \partial_x^4 U_0 \nonumber \\
    & +3 \tau_0^3 \alpha  \gamma^2 U_0 (\partial_x U_0)^2 \partial_x^4 U_0 +\frac{1}{2} \tau_0^3 \alpha \gamma^2 U_0^2 \partial_x^2 U_0 \partial_x^4 U_0 +\frac{3}{2} \tau_0^3 \alpha^2 \gamma^2 (\partial_x^2 U_0)^2 \partial_x^4 U_0 + 22 \tau_0^3 \alpha^2 \gamma^2 \partial_x U_0 \partial_x^3 U_0 \partial_x^4 U_0 \nonumber \\
    & + 18 \tau_0^3 \alpha^2 \gamma^2 U_0 (\partial_x^4 U_0)^2 + \frac{1}{2} \tau_0 \alpha U_0 \partial_x^4 U_2 + \tau_0 \alpha \gamma U_0 \partial_x^4 U_2 -\tau_0 \alpha \gamma^2 U_0 \partial_x^4 U_2 \nonumber \\
    & + \frac{1}{2} \tau_0^3 \alpha \gamma^2 U_0^2 \partial_x U_0 \partial_x^5 U_0 + \frac{21}{2} \tau_0^3 \alpha^2 \gamma^2 \partial_x U_0 \partial_x^2 U_0 \partial_x^5 U_0 +\frac{57}{2} \tau_0^3 \alpha^2 \gamma^2 U_0 \partial_x^3 U_0 \partial_x^5 U_0 +126 \tau_0^3 \alpha^3 \gamma^3 (\partial_x^5 U_0)^2 \nonumber \\ 
    &-126 \tau_0^3 \alpha^3 \gamma^4 (\partial_x^5 U_0)^2 + 4 \tau_0^3 \alpha^2 \gamma^2 (\partial_x U_0)^2 \partial_x^6 U_0 + \frac{29}{2} \tau_0^3 \alpha^2 \gamma^2 U_0 \partial_x^2 U_0 \partial_x^6 U_0 +\frac{1}{2} \tau_0^3 \alpha^3 \gamma^2 \partial_x^4 U_0 \partial_x^6 U_0 \nonumber \\
    &+ 210 \tau_0^3 \alpha^3 \gamma^3 \partial_x^4 U_0 \partial_x^6 U_0 - 210 \tau_0^3 \alpha^3 \gamma^4 \partial_x^4 U_0 \partial_x^6 U_0 + \tau_0 \alpha^2 \gamma \partial_x^6 U_2 -\tau_0 \alpha^2 \gamma^2 \partial_x^6 U_2\nonumber \\
    & + 5 \tau_0^3 \alpha^2 \gamma^2 U_0 \partial_x U_0 \partial_x^7 U_0 + \frac{1}{2} \tau_0^3 \alpha^3 \gamma^2 \partial_x^3 U_0 \partial_x^7 U_0 + 120 \tau_0^3 \alpha^3 \gamma^3 \partial_x^3 U_0 \partial_x^7 U_0 -120 \tau_0^3 \alpha^3 \gamma^4 \partial_x^3 U_0 \partial_x^7 U_0 \nonumber \\
    &+ \frac{1}{2} \tau_0^3 \alpha^2 \gamma^2 U_0^2 \partial_x^8 U_0 + 45 \tau_0^3 \alpha^3 \gamma^3 \partial_x^2 U_0 \partial_x^8 U_0 - 45 \tau_0^3 \alpha^3 \gamma^4 \partial_x^2 U_0 \partial_x^8 U_0 +\frac{1}{2} \tau_0^3 \alpha^3 \gamma^2 \partial_x U_0 \partial_x^9 U_0 \nonumber \\ 
    & + 10 \tau_0^3 \alpha^3 \gamma^3 \partial_x U_0 \partial_x^9 U_0 - 10 \tau_0^3 \alpha^3 \gamma^4 \partial_x U_0 \partial_x^9 U_0  + \frac{1}{2} \tau_0^3 \alpha^3 \gamma^2 U_0 \partial_x^{10} U_0 +\tau_0^3 \alpha^3 \gamma^3 U_0 \partial_x^{10} U_0 \nonumber \\
    &- \tau_0^3 \alpha^3 \gamma^4 U_0  \partial_x^{10} U_0 + \tau_0^3 \alpha^4 \gamma^3 \partial_x^{12} U_0 - \tau_0^3 \alpha^4 \gamma^4 \partial_x^{12} U_0, \nonumber
\end{align}
\end{subequations}
The primary (i.e., non-internal stage) equation given at each order invokes the following: 
\begin{itemize}
    \item \eqref{eq:rk222-i2-e2} adopts $\partial_{t_1} U_0 = 0$ and $U_1 = 0$
    \item \eqref{eq:rk222-i2-e3} adopts $\partial_{t_1} U_0 = 0$, $U_1 = 0$, $\partial_{t_2} U_0 = 0$, and $\partial_{t_1} U_2 = 0$
\end{itemize}
For the internal stage equations, however, we do not make these assumptions explicit, though we could do so without changing the end result. We obtain the nontrivial solvability condition by integrating \eqref{eq:rk222-i2-e3} with respect to $\xi = x - ct_0$.

\subsection{RK443}

For RK443, the multiple scales analysis yields
\begin{subequations}
\begin{align}
    {\cal O}(\epsilon^{-1}) \quad &U_0^{(1)} = U_0, \label{eq:rk443-i1-en1} \\
    &U_0^{(2)} = U_0, \label{eq:rk443-i2-en1}\\
    &U_0^{(3)} = U_0, \label{eq:rk443-i3-en1}\\
    {\cal O}(\epsilon^0) \quad &\partial_{t_0} U_0 + \alpha \partial_x^3 U_0 + U_0 \partial_x U_0 = 0, \label{eq:rk443-i4-e0} \\
    &U_1^{(1)} = U_1 - \frac{1}{2} \tau_0 U_0 \partial_x U_0 - \frac{1}{2} \tau_0 \alpha \partial_x^3 U_0^{(1)}, \label{eq:rk443-i1-e0}\\
    & U_1^{(2)} = U_1 - \frac{1}{18}\tau_0 U_0^{(1)} \partial_x U_0^{(1)} - \frac{11}{18} \tau_0 U_0 \partial_x U_0 - \frac{1}{6} \tau_0 \alpha \partial_x^3 U_0^{(1)} - \frac{1}{2} \tau_0 \alpha \partial_x^3 U_0^{(2)}, \label{eq:rk443-i2-e0}\\
    & U_1^{(3)} = U_1 + \frac{5}{6} \tau_0 U_0^{(1)} \partial_x U_0^{(1)} - \frac{1}{2} \tau_0 U_0^{(2)} \partial U_0^{(2)} - \frac{5}{6} \tau_0 U_0 \partial_x U_0 + \frac{1}{2} \tau_0 \alpha \partial_x^3 U_0^{(1)} - \frac{1}{2} \tau_0 \alpha \partial_x^3 U_0^{(2)} - \frac{1}{2} \tau_0 \alpha \partial_x^3 U_0^{(3)}, \label{eq:rk443-i3-e0}\\
    {\cal O}(\epsilon^1) \quad &\partial_{t_0}U_1 + \alpha \partial_x^3 U_1 + \partial_x (U_0 U_1) = -\partial_{t_1} U_0, \label{eq:rk443-i4-e1}\\
    &U_2^{(1)} = U_2 - \frac{1}{2} \tau_0 U_1 \partial_x U_0 - \frac{1}{2} \tau_0 U_0 \partial_x U_1 - \frac{1}{2} \tau_0 \alpha \partial_x^3 U_1^{(1)}, \label{eq:rk443-i1-e1}\\
    &U_2^{(2)} = U_2 - \frac{1}{18} \tau_0 U_1^{(1)} \partial_x U_0^{(1)} - \frac{1}{18} \tau_0 U_0^{(1)} \partial_x U_1^{(1)} - \frac{11}{18} \tau_0 U_1 \partial_x U_0 - \frac{11}{18} \tau_0 U_0 \partial_x U_1 - \frac{1}{6} \tau_0 \alpha \partial_x^3 U_1^{(1)} - \frac{1}{2} \tau_0 \alpha \partial_x^3 U_1^{(2)}, \label{eq:rk443-i2-e1}\\
    &U_2^{(3)} =  U_2 + \frac{5}{6} \tau_0 U_1^{(1)} \partial_x U_0^{(1)} + \frac{5}{6} \tau_0 U_0^{(1)} \partial_x U_1^{(1)} - \frac{1}{2} \tau_0 U_1^{(2)} \partial_X U_0^{(2)} - \frac{1}{2} \tau_0 U_0^{(2)} \partial_x U_1^{(2)} - \frac{5}{6} \tau_0 U_1 \partial_x U_0 \label{eq:rk443-i3-e1} \\
    &- \frac{5}{6} \tau_0 U_0 \partial_x U_1 + \frac{1}{2} \tau_0 \alpha \partial_x^3 U_1^{(1)} - \frac{1}{2} \tau_0 \alpha \partial_x^3 U_1^{(2)} - \frac{1}{2} \tau_0 \alpha \partial_x^3 U_1^{(3)}, \nonumber \\
    {\cal O}(\epsilon^2) \quad &\partial_{t_0} U_2 + \alpha \partial_x^3 U_2 + \partial_x(U_0 U_2) = - \partial_{t_2} U_0, \label{eq:rk443-i4-e2}\\
    &U_3^{(1)} =  U_3 - \frac{1}{2}\tau_0 U_2 \partial_x U_0 -\frac{1}{2} \tau_0 U_1 \partial_x U_1 
    - \frac{1}{2}\tau_0 U_0 \partial_x U_2 - \frac{1}{2} \tau_0 \alpha \partial_x^3 U_2^{(1)}, \label{eq:rk443-i1-e2}\\
    &U_3^{(2)} =  U_3 - \frac{1}{18} \tau_0 U_2^{(1)}\partial_x U_0^{(1)} - \frac{1}{18} \tau_0 U_1^{(1} \partial_x U_1^{(1)} - \frac{1}{18} \tau_0 U_0^{(1)} \partial_x U_2^{(1)} - \frac{11}{18} \tau_0 U_2 \partial_x U_0 - \frac{11}{18} \tau_0 U_1 \partial_x U_1 \label{eq:rk443-i2-e2}\\
    &- \frac{11}{18} \tau_0 U_0 \partial_x U_2 - \frac{1}{6} \tau_0 \alpha \partial_x^3 U_2^{(1)} - \frac{1}{2} \tau_0 \alpha  \partial_x^3 U_2^{(2)}, \nonumber \\
    &U_3^{(3)} = U_3 + \frac{5}{6} \tau_0 U_2^{(1)} \partial_x U_0^{(1)} + \frac{5}{6} \tau_0 U_1^{(1)} \partial_x U_1^{(1)} + \frac{5}{6} \tau_0 U_0^{(1)} \partial_x U_2^{(1)} - \frac{1}{2} \tau_0 U_2^{(2)} \partial_x U_0^{(2)} - \frac{1}{2} \tau_0 U_1^{(2)} \partial_x U_1^{(2)} \label{eq:rk443-i3-e2} \\
    &- \frac{1}{2} \tau_0 U_0^{(2)} \partial_x U_2^{(2)} - \frac{5}{6} \tau_0 U_2 \partial_x U_0 - \frac{5}{6} \tau_0 U_1 \partial_x U_1 - \frac{5}{6} \tau_0 U_0 \partial_x U_2 + \frac{1}{2} \tau_0 \alpha \partial_x^3 U_2^{(1)} - \frac{1}{2} \tau_0 \alpha \partial_x^3 U_2^{(2)} - \frac{1}{2} \tau_0 \alpha \partial_x^3 U_2^{(3)},\nonumber \\
    {\cal O}(\epsilon^3) \quad &\partial_{t_0}U_3 + \alpha \partial_x^3 U_3 + \partial_x (U_0 U_3) = - \partial_{t_3}U_0 -\frac{1}{24} \tau_0^3 \partial_{t_0}^4 U_0 +\frac{1}{8} \tau_0 ^3 U_0 (\partial_x U_0)^4 - \frac{5}{48} \tau_0^3 U_0^2 (\partial_x U_0)^2 \partial_x^2 U_0 \label{eq:rk443-i4-e3}\\
    &-\frac{7}{72} \tau_0^3 U_0^3 (\partial_x^2 U_0)^2 + 3 \tau_0^3 \alpha (\partial_x U_0)^2 (\partial_x^2 U_0)^2 + \frac{5}{2} \tau_0^3 \alpha U_0 (\partial_x^2 U_0)^3 - \frac{7}{36} \tau_0^3 U_0^3 \partial_x U_0 \partial_x^3 U_0 \nonumber \\
    & + \frac{19}{8} \tau_0^3 \alpha  (\partial_x U_0)^3 \partial_x^3 U_0 + \frac{37}{3} \tau_0^3 \alpha U_0 \partial_x U_0 \partial_x^2 U_0 \partial_x^3 U_0 + \frac{35}{24} \tau_0^3 \alpha U_0^2 (\partial_x^3 U_0)^2 - \frac{645}{32} \tau_0^3 \alpha^2 \partial_x^2 U_0 (\partial_x^3 U_0)^2 \nonumber \\
    &-\frac{1}{12} \tau_0^3 \alpha \partial_{t_0}^3 \partial_x^3 U_0 -\frac{7}{288} \tau_0^3 U_0^4 \partial_x^4 U_0 +\frac{63}{16} \tau_0^3 \alpha U_0 (\partial_x U_0)^2 \partial_x^4 U_0 + \frac{223}{96} \tau_0^3 \alpha U_0^2 \partial_x^2 U_0 \partial_x^4 U_0 \nonumber \\
    &-\frac{517}{32} \tau_0^3 \alpha^2 (\partial_x^2 U_0)^2 \partial_x^4 U_0 -\frac{3209}{144} \tau_0^3 \alpha^2 \partial_x U_0 \partial_x^3 U_0 \partial_x^4 U_0 - \frac{545}{144} \tau_0^3 \alpha^2 U_0 (\partial_x^4 U_0)^2 + \frac{49}{48} \tau_0^3 \alpha U_0^2 \partial_x U_0 \partial_x^5 U_0 \nonumber \\
    &-\frac{1385}{96} \tau_0^3 \alpha^2 \partial_x U_0 \partial_x^2 U_0 \partial_x^5 U_0 -\frac{1735}{288} \tau_0^3 \alpha^2 U_0 \partial_x^3 U_0 \partial_x^5 U_0 -\frac{415}{48} \tau_0^3 \alpha^3 (\partial_x^5 U_0)^2 + \frac{5}{96} \tau_0^3 \alpha U_0^3 \partial_x^6 U_0 \nonumber \\
    &-\frac{665}{288} \tau_0^3 \alpha^2 (\partial_x U_0)^2 \partial_x^6 U_0 -\frac{197}{72} \tau_0^3 \alpha^2 U_0 \partial_x^2 U_0 \partial_x^6 U_0 -\frac{4183}{288} \tau_0^3 \alpha^3 \partial_x^4 U_0 \partial_x^6 U_0 -\frac{149}{288} \tau_0^3 \alpha^2 U_0 \partial_x U_0 \partial_x^7 U_0 \nonumber \\
    &-\frac{2509}{288} \tau_0^3 \alpha^3 \partial_x^3 U_0 \partial_x^7 U_0 -\frac{1}{96} \tau_0^3 \alpha^2 U_0^2 \partial_x^8 U_0 -\frac{11}{3} \tau_0^3 \alpha^3 \partial_x^2 U_0 \partial_x^8 U_0 -\frac{283}{288} \tau_0^3 \alpha^3 \partial_x U_0 \partial_x^9 U_0 \nonumber \\
    &-\frac{43}{288} \tau_0^3 \alpha^3 U_0 \partial_x^{10} U_0 -\frac{1}{16} \tau_0^3 \alpha^4 \partial_x^{12} U_0\nonumber.
\end{align}
\end{subequations}
The primary (i.e., non-internal stage) equation given at each order invokes the following: 
\begin{itemize}
    \item \eqref{eq:rk443-i4-e2} adopts $\partial_{t_1} U_0 = 0$ and $U_1 = 0$
    \item \eqref{eq:rk443-i4-e3} adopts $\partial_{t_1} U_0 = 0$, $U_1 = 0$, $\partial_{t_2} U_0 = 0$, and $U_2 = 0$
\end{itemize}
As described for RK222, for the internal stage equations we do not make these assumptions explicit, though we could do so without changing the end result. We obtain the nontrivial solvability condition by integrating \eqref{eq:rk443-i4-e3} with respect to $\xi = x - ct_0$. 

\section{Evolution of the KdV soliton parameter $c$ from nontrivial solvability conditions}\label{sec:sc}

We derive particular expressions of nontrivial solvability conditions for the case that the initial data of the simulation is given by a soliton of the form
\begin{equation}\label{eq:soliton-sm}
    u(x,0) = u_{\infty}(c_0) + 3c_0\,{\rm sech}\left(\sqrt{\frac{c_0}{4\alpha}}x\right)^2, 
\end{equation} 
where $u_{\infty}$ is defined as
\begin{equation}\label{eq:ms-uinf-sm}
    u_{\infty} \equiv u_{\infty}(c) =\frac{12\alpha^{1/2}}{L}\left(c_0^{1/2}{\rm tanh}{\left(\frac{L}{2}\sqrt{\frac{c_0}{4\alpha}}\right)} - c^{1/2}{\rm tanh}{\left(\frac{L}{2}\sqrt{\frac{c}{4\alpha}}\right)}\right). 
\end{equation}
While $u_{\infty}(c_0)$ is identically zero, we need the expression~\eqref{eq:ms-uinf-sm} for a general value of $c$, as $c$ will eventually be found to vary with respect to a slow time scale. 

({\textbf NOTE:} to make our numerical integration of solvability conditions more efficient, we choose to work with the approximation of~\eqref{eq:ms-uinf-sm} that
\begin{equation}\label{eq:ms-uinf-approx-sm}
    u_{\infty}(c) \approx \frac{12\alpha^{1/2}}{L}\left(c_0^{1/2} - c^{1/2}\right),
\end{equation}
which is valid in the regime of large $L$ and small $\alpha$ we work in. For SBDF1, SBDF2, and RK222, $c$ is furthermore an increasing function of time, so the approximation remains valid as $t \rightarrow \infty$. This is not the case for RK443 where $c$ decreases in time, but the approximation is still highly accurate provided that $c$ does not get too small. As our RK443 simulations never go beyond the point that $c$ has decreased by $\approx$ 20\% of its original value---which takes an extremely long time to achieve---the approximation~\eqref{eq:ms-uinf-approx-sm} is still very accurate, with the discrepancy between~\eqref{eq:ms-uinf-sm} and~\eqref{eq:ms-uinf-approx-sm} being several orders of magnitude smaller than the smallest errors being plotted in the numerical results of the paper.)

In the $L \rightarrow \infty$ limit, the respective nontrivial solvability conditions for each scheme are given by the separable ODEs
\begin{equation}\label{eq:scs-infinite}
\begin{aligned}
    \frac{\partial c}{\partial t_1} &= \frac{34}{105}\frac{\tau_0}{\alpha}c^4, \quad \text{(SBDF1)}\\
    \frac{\partial c}{\partial t_3} &= \frac{43}{105}\frac{\tau_0^3}{\alpha^2}c^7, \quad \text{(SBDF2)}\\
    \frac{\partial c}{\partial t_3} &= \frac{355045-245436\sqrt{2}}{30030}\frac{\tau_0^3}{\alpha^2}c^7, \quad \text{(RK222)}\\
    \frac{\partial c}{\partial t_3} &= -\frac{77069}{180180}\frac{\tau_0^3}{\alpha^2}c^7. \quad \text{(RK443)}\\
\end{aligned}
\end{equation}
When we work on the periodic domain with finite $L$, the nontrivial solvability conditions result in ODEs of the form 
\begin{equation}\label{eq:sc-finite}
    \frac{\partial c}{\partial t_{m^*}} = \frac{{\rm RHS}_{\rm scheme}(c; \alpha, L, c_0)}{{\rm LHS}(c;\alpha, L, c_0)},
\end{equation}
where $m^* = 1$ for SBDF1 and $m^* = 3$ for SBDF2, RK222 and RK443. In general we would expect that $m^*$ would be given by the order of accuracy of the scheme, but (as discussed in the main paper) parity considerations in the context of KdV solitons cause the nontrivial conditions of the second order schemes to first appear at higher orders in $\epsilon$. The expressions for ${\rm LHS}$ and ${\rm RHS}_{\rm scheme}$ are given below.
\begin{align}\label{eq:LHS}
    {\rm LHS}(c;\alpha, L, c_0) = \frac{9}{4 L \sqrt{c}} \Biggl[&L^2 c^{3/2} {\rm sech}\left(\frac{L \sqrt{c}}{4 \sqrt{\alpha}}\right)^4+32 \sqrt{c_0} \alpha\left(-1+{\rm tanh}\left(\frac{L \sqrt{c}}{4\sqrt{\alpha}}\right)\right)\\
    &+8 \sqrt{c} \left(4 \alpha + \sqrt{c_0} L \sqrt{\alpha } {\rm sech}\left(\frac{L\sqrt{c}}{4 \sqrt{\alpha }}\right)^2-8 \alpha  {\rm tanh}\left(\frac{L \sqrt{c}}{4 \sqrt{\alpha}}\right)\right) \nonumber \\
    &-4 L \sqrt{\alpha} c \left(2{\rm sech}\left(\frac{L \sqrt{c}}{4 \sqrt{\alpha }}\right)^2-3 {\rm tanh}\left(\frac{L \sqrt{c}}{4 \sqrt{\alpha}}\right)+{\rm tanh}\left(\frac{L\sqrt{c}}{4\sqrt{\alpha}}\right)^3\right)\Biggr], \nonumber%\frac{\partial c}{\partial t_1}
\end{align}

\begin{align}\label{eq:RHS-SBDF1}
    {\rm RHS}_{\rm SBDF1}(c; \alpha, L, c_0) = &\frac{3 c^{3/2} \tau_0 \left(c L+12 \left(-\sqrt{c}+\sqrt{c_0}\right) \sqrt{\alpha }\right)}{35 L^3 \sqrt{\alpha }} \Biggl[4 c L \left(17c L+84 \left(-\sqrt{c}+\sqrt{c_0}\right) \sqrt{\alpha }\right)  \\
    &+2 \biggl(17 c^2 L^2+126 c \left(\sqrt{c}-\sqrt{c_0}\right) L \sqrt{\alpha} +2520 \left(c-2 \sqrt{c} \sqrt{c_0}+c_0\right) \alpha \biggr) {\rm sech}\left(\frac{\sqrt{c} L}{4 \sqrt{\alpha }}\right)^2 \nonumber \\
    &-9 c L \left(3c L+364 \left(\sqrt{c}-\sqrt{c_0}\right) \sqrt{\alpha }\right) {\rm sech}\left(\frac{\sqrt{c} L}{4 \sqrt{\alpha }}\right)^4 +450 c^2 L^2 {\rm sech}\left(\frac{\sqrt{c}L}{4 \sqrt{\alpha }}\right)^6\Biggr] {\rm tanh}\left(\frac{\sqrt{c} L}{4 \sqrt{\alpha }}\right),\nonumber
\end{align}

\begin{align}\label{eq:RHS-SBDF2}
    {\rm RHS}_{\rm SBDF2}(c; \alpha, L, c_0) = &\frac{3 \tau_0^3 c^{5/2} \left(12 \sqrt{c_0} \sqrt{\alpha }-12 \sqrt{\alpha } \sqrt{c}+L c\right)^3}{70 L^5 \alpha ^{3/2}} \Biggl[720 \sqrt{c_0} L \sqrt{\alpha } c-720 L \sqrt{\alpha } c^{3/2}+172 L^2 c^2\\
    &-15120 c_0 \alpha  {\rm sech}\left(\frac{L\sqrt{c}}{4 \sqrt{\alpha }}\right)^2+30240 \sqrt{c_0} \alpha  \sqrt{c} {\rm sech}\left(\frac{L \sqrt{c}}{4 \sqrt{\alpha }}\right)^2+780\sqrt{c_0} L \sqrt{\alpha } c {\rm sech}\left(\frac{L \sqrt{c}}{4 \sqrt{\alpha }}\right)^2\nonumber \\
    &-15120 \alpha  c {\rm sech}\left(\frac{L\sqrt{c}}{4 \sqrt{\alpha }}\right)^2-780 L \sqrt{\alpha } c^{3/2} {\rm sech}\left(\frac{L \sqrt{c}}{4 \sqrt{\alpha }}\right)^2+86 L^2 c^2{\rm sech}\left(\frac{L \sqrt{c}}{4 \sqrt{\alpha }}\right)^2\nonumber \\
    &+45360 c_0 \alpha  {\rm sech}\left(\frac{L \sqrt{c}}{4 \sqrt{\alpha }}\right)^4-90720\sqrt{c_0} \alpha  \sqrt{c} {\rm sech}\left(\frac{L \sqrt{c}}{4 \sqrt{\alpha }}\right)^4-23040 \sqrt{c_0} L \sqrt{\alpha } c{\rm sech}\left(\frac{L \sqrt{c}}{4 \sqrt{\alpha }}\right)^4\nonumber \\
    &+45360 \alpha  c {\rm sech}\left(\frac{L \sqrt{c}}{4 \sqrt{\alpha }}\right)^4+23040L \sqrt{\alpha } c^{3/2} {\rm sech}\left(\frac{L \sqrt{c}}{4 \sqrt{\alpha }}\right)^4+117 L^2 c^2 {\rm sech}\left(\frac{L \sqrt{c}}{4\sqrt{\alpha }}\right)^4\nonumber \\
    &+45900 \sqrt{c_0} L \sqrt{\alpha } c {\rm sech}\left(\frac{L \sqrt{c}}{4 \sqrt{\alpha }}\right)^6-45900 L \sqrt{\alpha} c^{3/2} {\rm sech}\left(\frac{L \sqrt{c}}{4 \sqrt{\alpha }}\right)^6-3525 L^2 c^2 {\rm sech}\left(\frac{L \sqrt{c}}{4 \sqrt{\alpha}}\right)^6\nonumber \\
    &+7350 L^2 c^2 {\rm sech}\left(\frac{L \sqrt{c}}{4 \sqrt{\alpha }}\right)^8\Biggr] {\rm tanh}\left(\frac{L \sqrt{c}}{4 \sqrt{\alpha}}\right),\nonumber
\end{align}

\begin{align}\label{eq:RHS-RK222}
   {\rm RHS}_{\rm RK222}(c; \alpha, L, c_0) = &\frac{3 \tau_0^3 c^{5/2} \left(12 \sqrt{c_0} \sqrt{\alpha }-12 \sqrt{\alpha } \sqrt{c}+L c\right)}{40040 L^5 \alpha^{3/2}}\Biggl[-8 L c \Biggl(-1235520 c_0^{3/2} \alpha^{3/2}+3706560 c_0 \alpha^{3/2} \sqrt{c}\\
   &+20592 \left(\left(-71+28 \sqrt{2}\right)c_0 L \alpha -180 \sqrt{c_0} \alpha ^{3/2}\right) c-41184 \left(\left(-71+28 \sqrt{2}\right) \sqrt{c_0} L \alpha -30 \alpha^{3/2}\right) c^{3/2}\nonumber \\
   &+156 \left(\left(-6599+4116 \sqrt{2}\right) \sqrt{c_0} L^2 \sqrt{\alpha }+132 \left(-71+28 \sqrt{2}\right) L \alpha\right) c^2\nonumber \\
   &-156 \left(-6599+4116 \sqrt{2}\right) L^2 \sqrt{\alpha } c^{5/2}+\left(-355045+245436 \sqrt{2}\right) L^3 c^3\Biggr)\nonumber \\
   &-4 \Biggl(51891840c_0^2 \alpha ^2-207567360 c_0^{3/2} \alpha^2 \sqrt{c}+617760 \left(\left(-37+28 \sqrt{2}\right) c_0^{3/2} L \alpha^{3/2}+504c_0 \alpha^2\right) c\nonumber \\
   &-1853280 \left(\left(-37+28 \sqrt{2}\right) c_0 L \alpha^{3/2}+112 \sqrt{c_0} \alpha^2\right) c^{3/2}\nonumber \\
   &-10296\left(\left(-453+364 \sqrt{2}\right) c_0 L^2-180 \left(-37+28 \sqrt{2}\right) \sqrt{c_0} L \sqrt{\alpha }-5040 \alpha \right) \alpha c^2\nonumber \\
   &+20592 L \left(\left(-453+364 \sqrt{2}\right) \sqrt{c_0} L+30 \left(37-28 \sqrt{2}\right) \sqrt{\alpha }\right) \alpha  c^{5/2}\nonumber \\
   &+234\left(\left(-6581+4284 \sqrt{2}\right) \sqrt{c_0} L^3 \sqrt{\alpha }+44 \left(453-364 \sqrt{2}\right) L^2 \alpha \right) c^3\nonumber \\
   &-234 \left(-6581+4284\sqrt{2}\right) L^3 \sqrt{\alpha } c^{7/2}+\left(-355045+245436 \sqrt{2}\right) L^4 c^4\Biggr) {\rm sech}\left(\frac{L \sqrt{c}}{4 \sqrt{\alpha}}\right)^2\nonumber \\
   &-18 \Biggl(-34594560 c_0^2 \alpha^2+138378240 c_0^{3/2} \alpha ^2 \sqrt{c}-823680 \left(\left(-110+63 \sqrt{2}\right) c_0^{3/2}L \alpha^{3/2}+252 c_0 \alpha ^2\right) c\nonumber \\
   &+2471040 \left(\left(-110+63 \sqrt{2}\right) c_0 L \alpha^{3/2}+56 \sqrt{c_0} \alpha^2\right) c^{3/2}\nonumber \\
   &+6864 \left(\left(-14998+10899 \sqrt{2}\right) c_0 L^2-360 \left(-110+63 \sqrt{2}\right) \sqrt{c_0} L \sqrt{\alpha}-5040 \alpha \right) \alpha  c^2\nonumber \\
   &-13728 L \left(\left(-14998+10899 \sqrt{2}\right) \sqrt{c_0} L+60 \left(110-63 \sqrt{2}\right) \sqrt{\alpha}\right) \alpha  c^{5/2}\nonumber \\
   &-52 \left(\left(-744948+526547 \sqrt{2}\right) \sqrt{c_0} L^3 \sqrt{\alpha }+132 \left(14998-10899 \sqrt{2}\right)L^2 \alpha \right) c^3 \nonumber \\
   &+52 \left(-744948+526547 \sqrt{2}\right) L^3 \sqrt{\alpha } c^{7/2}+\left(-73355+50916 \sqrt{2}\right) L^4 c^4\Biggr){\rm sech}\left(\frac{L \sqrt{c}}{4 \sqrt{\alpha }}\right)^4\nonumber \\
   &+30 L c \Biggl(-741312 \left(-124+63 \sqrt{2}\right) c_0^{3/2} \alpha ^{3/2}+2223936\left(-124+63 \sqrt{2}\right) c_0 \alpha ^{3/2} \sqrt{c}\nonumber \\
   &+10296 \left(\left(-32782+23373 \sqrt{2}\right) c_0 L \alpha -216 \left(-124+63\sqrt{2}\right) \sqrt{c_0} \alpha ^{3/2}\right) c\nonumber \\
   &-20592 \left(\left(-32782+23373 \sqrt{2}\right) \sqrt{c_0} L \alpha +36 \left(124-63\sqrt{2}\right) \alpha ^{3/2}\right) c^{3/2}\nonumber \\
   &-78 \left(\left(-4321172+3064593 \sqrt{2}\right) \sqrt{c_0} L^2 \sqrt{\alpha }+132 \left(32782-23373\sqrt{2}\right) L \alpha \right) c^2\nonumber \\
   &+78 \left(-4321172+3064593 \sqrt{2}\right) L^2 \sqrt{\alpha } c^{5/2}+\left(-3845701+2717282 \sqrt{2}\right)L^3 c^3\Biggr) {\rm sech}\left(\frac{L \sqrt{c}}{4 \sqrt{\alpha }}\right)^6\nonumber \\
   &-105 L^2 c^2 \Biggl(20592 \left(-5002+3501 \sqrt{2}\right) c_0\alpha -41184 \left(-5002+3501 \sqrt{2}\right) \sqrt{c_0} \alpha  \sqrt{c}\nonumber \\
   &-312 \left(\left(-950770+674649 \sqrt{2}\right) \sqrt{c_0}L \sqrt{\alpha }+66 \left(5002-3501 \sqrt{2}\right) \alpha \right) c\nonumber \\
   &+312 \left(-950770+674649 \sqrt{2}\right) L \sqrt{\alpha } c^{3/2}+\left(-17580062+12458161\sqrt{2}\right) L^2 c^2\Biggr) {\rm sech}\left(\frac{L \sqrt{c}}{4 \sqrt{\alpha }}\right)^8\nonumber \\
   &+2835 L^3 c^3 \Biggl(208 \left(41908-29715 \sqrt{2}\right)\sqrt{c_0} \sqrt{\alpha }+208 \left(-41908+29715 \sqrt{2}\right) \sqrt{\alpha } \sqrt{c}\nonumber \\
   &+\left(-2084828+1478415 \sqrt{2}\right) L c\Biggr){\rm sech}\left(\frac{L \sqrt{c}}{4 \sqrt{\alpha }}\right)^{10}\nonumber \\
   &-1715175 \left(-2788+1977 \sqrt{2}\right) L^4 c^4 {\rm sech}\left(\frac{L\sqrt{c}}{4 \sqrt{\alpha }}\right)^{12}\Biggr] {\rm tanh}\left(\frac{L \sqrt{c}}{4 \sqrt{\alpha }}\right),\nonumber
\end{align}

\begin{align}\label{eq:RHS-RK443}
   {\rm RHS}_{\rm RK443}(c; \alpha, L, c_0) = &-\frac{\tau_0^3 c^{5/2} \left(12 \sqrt{c_0} \sqrt{\alpha }-12 \sqrt{\alpha } \sqrt{c}+L c\right)}{160160 L^5 \alpha ^{3/2}}\Biggl[16 L c \Biggl(3912480 c_0^{3/2} \alpha ^{3/2}-11737440 c_0 \alpha ^{3/2} \sqrt{c}\\
   &+6864 \left(377 c_0 L \alpha +1710\sqrt{c_0} \alpha ^{3/2}\right) c-13728 \left(377 \sqrt{c_0} L \alpha +285 \alpha ^{3/2}\right) c^{3/2}\nonumber \\
   &+26 \left(27619 \sqrt{c_0}L^2 \sqrt{\alpha }+99528 L \alpha \right) c^2-718094 L^2 \sqrt{\alpha } c^{5/2}+77069 L^3 c^3\Biggr)\nonumber \\
   &+8 \Biggl(-164324160 c_0^2 \alpha^2+657296640 c_0^{3/2} \alpha ^2 \sqrt{c}+3706560 \left(3 c_0^{3/2} L \alpha ^{3/2}-266 c_0 \alpha ^2\right) c\nonumber \\
   &-1235520 \left(27c_0 L \alpha ^{3/2}-532 \sqrt{c_0} \alpha ^2\right) c^{3/2}+15444 \left(117 c_0 L^2+2160 \sqrt{c_0} L \sqrt{\alpha }-10640\alpha \right) \alpha  c^2\nonumber \\
   &-277992 L \left(13 \sqrt{c_0} L+40 \sqrt{\alpha }\right) \alpha  c^{5/2}+52 \left(13232 \sqrt{c_0} L^3\sqrt{\alpha }+34749 L^2 \alpha \right) c^3-688064 L^3 \sqrt{\alpha } c^{7/2}\nonumber \\
   &+77069 L^4 c^4\Biggr) {\rm sech}\left(\frac{L \sqrt{c}}{4\sqrt{\alpha }}\right)^2+24 \Biggl(164324160 c_0^2 \alpha ^2-657296640 c_0^{3/2} \alpha ^2 \sqrt{c}\nonumber \\
   &-926640 \left(179 c_0^{3/2}L \alpha ^{3/2}-1064 c_0 \alpha ^2\right) c+308880 \left(1611 c_0 L \alpha ^{3/2}-2128 \sqrt{c_0} \alpha ^2\right) c^{3/2}\nonumber \\
   &+5148\alpha  \left(8479 c_0 L^2-96660 \sqrt{c_0} L \sqrt{\alpha }+31920 \alpha \right) c^2-10296 L \left(8479 \sqrt{c_0} L-16110\sqrt{\alpha }\right) \alpha  c^{5/2}\nonumber \\
   &+13 \left(299287 \sqrt{c_0} L^3 \sqrt{\alpha }+3357684 L^2 \alpha \right) c^3-3890731 L^3 \sqrt{\alpha} c^{7/2}+18016 L^4 c^4\Biggr) {\rm sech}\left(\frac{L \sqrt{c}}{4 \sqrt{\alpha }}\right)^4\nonumber \\
   &+30 L c \Biggl(264771936 c_0^{3/2} \alpha^{3/2}-794315808 c_0 \alpha ^{3/2} \sqrt{c}-61776 \left(4189 c_0 L \alpha -12858 \sqrt{c_0} \alpha ^{3/2}\right) c\nonumber \\
   &+123552\left(4189 \sqrt{c_0} L \alpha -2143 \alpha ^{3/2}\right) c^{3/2}-78 \left(270913 \sqrt{c_0} L^2 \sqrt{\alpha }+3317688 L \alpha \right)c^2\nonumber \\
   &+21131214 L^2 \sqrt{\alpha } c^{5/2}+512177 L^3 c^3\Biggr) {\rm sech}\left(\frac{L \sqrt{c}}{4 \sqrt{\alpha }}\right)^6-105 L^2 c^2\Biggl(-92849328 c_0 \alpha\nonumber \\
   &+185698656 \sqrt{c_0} \alpha  \sqrt{c}-312 \left(1715 \sqrt{c_0} L \sqrt{\alpha }+297594 \alpha \right)c+535080 L \sqrt{\alpha } c^{3/2}\nonumber \\
   &+ 1449353 L^2 c^2\Biggr) {\rm sech}\left(\frac{L \sqrt{c}}{4 \sqrt{\alpha }}\right)^8+1890 L^3 c^3\Biggl(622492 \sqrt{c_0} \sqrt{\alpha }-622492 \sqrt{\alpha } \sqrt{c}\nonumber \\
   &+150691 L c\Biggr) {\rm sech}\left(\frac{L \sqrt{c}}{4 \sqrt{\alpha}}\right)^{10}-112058100 L^4 c^4 {\rm sech}\left(\frac{L \sqrt{c}}{4 \sqrt{\alpha }}\right)^{12}\Biggr] {\rm tanh}\left(\frac{L \sqrt{c}}{4\sqrt{\alpha}}\right).\nonumber
\end{align}

\section{Extension of the multiple scales analysis to fully-implicit schemes}\label{sec:implicit}

Here we briefly consider how the multiple scales analysis can be applied to some fully implicit timestepping schemes without any significant modification of the approach. We consider the fully-implicit analogs of SBDF1 and SBDF2, i.e., BDF1 (Backward-Euler) and BDF2. Written for KdV, BDF1 corresponds to the time discretization
\begin{equation}\label{eq:bdf1-kdv}
    \frac{1}{\Delta t}\left(U^{(n+1)}-U^{(n)}\right) + \alpha \partial_x^3 U^{(n+1)} + U^{(n+1)} \partial_x U^{(n+1)} = 0,
\end{equation}
and BDF2 similarly corresponds to
\begin{equation}\label{eq:bdf2-kdv}
    \frac{1}{\Delta t}\left(\frac{3}{2}U^{(n+2)}-2U^{(n+1)}+\frac{1}{2}U^{(n)}\right) + \alpha \partial_x^3 U^{(n+2)} + U^{(n+2)} \partial_x U^{(n+2)} = 0.
\end{equation}
Following the same setup and procedure described in~\ref{sec:ms-each}, here we report the predictions of the multiple scales analysis for these implicit schemes at each order in $\epsilon$.

For BDF1, the multiple scales analysis yields
\begin{subequations}
\begin{align}
    {\cal O}(\epsilon^0) \quad &\partial_{t_0} U_0 + \alpha \partial_x^3 U_0 + U_0 \partial_x U_0 = 0, \label{eq:bdf1-e0} \\
    {\cal O}(\epsilon^1) \quad &\partial_{t_0} U_1 + \alpha\partial_x^3 U_1 + \partial_x(U_0 U_1) = - \partial_{t_1} U_0 + \tau_0\partial_{t_0}\left(-\frac{\alpha}{2}\partial_x^3U_0 - \frac{1}{2}U_0\partial_xU_0\right). \label{eq:bdf1-e1}
\end{align}
\end{subequations}

For BDF2, we have
\begin{subequations}
\begin{align}
    {\cal O}(\epsilon^0) \quad &\partial_{t_0} U_0 + \alpha \partial_x^3 U_0 + U_0 \partial_x U_0 = 0, \label{eq:bdf2-e0} \\
    {\cal O}(\epsilon^1) \quad &\partial_{t_0} U_1 + \alpha \partial_x^3 U_1 + \partial_x \left(U_0 U_1\right) = -\partial_{t_1} U_0, \label{eq:bdf2-e1}\\
    {\cal O}(\epsilon^2) \quad &\partial_{t_0} U_2 + \alpha \partial_x^3 U_2 + \partial_x \left(U_0 U_2\right) = -\partial_{t_2} U_0 + \tau_0^2 \partial_{t_0}^2\left(-\frac{\alpha}{3}\partial_x^3 U_0 - \frac{1}{3}U_0\partial_x U_0\right), \label{eq:bdf2-e2} \\
    {\cal O}(\epsilon^3) \quad &\partial_{t_0}U_3 + \alpha \partial_x^3 U_3 + \partial_x (U_0 U_3) = - \partial_{t_3}U_0 + \tau_0^3 \partial_{t_0}^3 \left(\frac{\alpha}{4} \partial_x^3 U_0 + \frac{1}{4} U_0 \partial_x U_0 \right). \label{eq:bdf2-e3}
\end{align}
\end{subequations}
For each equation after \eqref{eq:bdf2-e1}, we adopt $\partial_{t_1} U_0 = 0$ and $U_1 = 0$. After $\eqref{eq:bdf2-e2}$ we further adopt $\partial_{t_2} U_0 = 0$ and $\partial_{t_1} U_2 = 0$.

Notably, the analysis for BDF1 and BDF2 nearly exactly mirrors that for SBDF1 and SBDF2. They differ in how the $\alpha\partial_x^3 U_0$ and $U_0\partial_x U_0$ terms in equations~\eqref{eq:bdf1-e1},~\eqref{eq:bdf2-e2}, and~\eqref{eq:bdf2-e3} now share the same coefficients, whereas they differ from one another in their counterpart equations~\eqref{eq:sbdf1-e1},~\eqref{eq:sbdf2-e2}, and~\eqref{eq:sbdf2-e3}.

Finally, following the same procedure as described in Section~\ref{sec:sc}, in the $L\rightarrow \infty$ limit, BDF1 and BDF2 admit the respective nontrivial solvability conditions given by the separable ODEs
\begin{equation}\label{eq:scs-infinite-im}
\begin{aligned}
    \frac{\partial c}{\partial t_1} &= -\frac{2}{15}\frac{\tau_0}{\alpha}c^4, \quad \text{(BDF1)}\\
    \frac{\partial c}{\partial t_3} &= -\frac{1}{21}\frac{\tau_0^3}{\alpha^2}c^7. \quad \text{(BDF2)}\\
\end{aligned}
\end{equation}

\printbibliography